\documentclass[12pt,leqno]{amsart}
\usepackage{amsmath, amssymb}
\usepackage{graphicx,color,hyperref}
\usepackage{verbatim, enumitem}

\newcommand\cA{{\mathcal A}}
\newcommand\cF{{\mathcal F}}
\newcommand\cG{{\mathcal G}}
\newcommand\cH{{\mathcal H}}
\newcommand\cI{{\mathcal I}}
\newcommand\cJ{{\mathcal J}}
\newcommand\cL{{\mathcal L}}
\newcommand\cN{{\mathcal N}}
\newcommand\cO{{\mathcal O}}

\newcommand\cS{{\mathcal S}}
\newcommand\cV{{\mathcal V}}
\newcommand\cW{{\mathcal W}}
\newcommand\cX{{\mathcal X}}
\newcommand\cY{{\mathcal Y}}
\newcommand\cZ{{\mathcal Z}}

\newcommand\bN{{\mathbb N}}
\newcommand\bC{{\mathbb C}}
\newcommand\bK{{\mathbb K}}
\newcommand\bZ{{\mathbb Z}}
\newcommand\bP{{\mathbb P}}
\newcommand\bQ{{\mathbb Q}}

\def\ii{{\rm i}\kern1pt}

\def\Ker{\mathop{\rm Ker}}

\def\Gr{\mathop{\rm Gr}\nolimits}
\def\GG{\mathop{\rm GG}\nolimits}
\def\GL{\mathop{\rm GL}\nolimits}

\def\Sing{\mathop{\rm Sing}}

\def\mod{\mathop{\rm mod}}
\def\rank{\mathop{\rm rank}}

\def\ECL{\mathop{\rm ECL}\nolimits}
\def\IEL{\mathop{\rm IEL}\nolimits}
\def\pr{\mathop{\rm pr}\nolimits}
\def\Pic{\mathop{\rm Pic}\nolimits}
\def\rank{\mathop{\rm rank}\nolimits}

\def\Sym{\mathop{\rm Sym}\nolimits}

\def\proof{\noindent{\it Proof.} }

\def\build#1^#2_#3{\mathrel{\mathop{\null#1}\limits^{#2}_{#3}}}
\def\mertorelbar{\vrule width0.6ex height0.65ex depth-0.55ex}
\def\merto{\mathrel{\mertorelbar\kern1.3pt\mertorelbar\kern1.3pt\mertorelbar
    \kern1.3pt\mertorelbar\kern-1ex\raise0.28ex\hbox{${\scriptscriptstyle>}$}}}

\catcode`\@=11
\newdimen\@rrowlength \@rrowlength=6ex
\def\ssrelbar{\vrule width\@rrowlength height0.64ex depth-0.56ex\kern-4pt}
\def\llra#1{\@rrowlength=#1\ssrelbar\rightarrow}
\catcode`\@=12

\def\semidirect{\mathop{\kern2pt\vrule depth-0.3pt height4.3pt 
\kern-2pt\times}\nolimits}

\def\HOMEPAGE{http://www-fourier.ujf-grenoble.fr/\~{}demailly/manuscripts}

\newdimen\plainitemindent \plainitemindent=18pt
\def\plainitem#1{\par\noindent
\hangindent\plainitemindent\hbox to\plainitemindent{#1\hss}\ignorespaces}

\catcode`\@=11
\def\openup{\afterassignment\@penup\dimen@=}
\def\@penup{\advance\lineskip\dimen@
  \advance\baselineskip\dimen@
  \advance\lineskiplimit\dimen@}
\newdimen\jot \jot=3pt
\newskip\plaincentering \plaincentering=0pt plus 1000pt minus 1000pt
\def\ialign{\everycr{}\tabskip\z@skip\halign}
\def\eqalign#1{\null\,\vcenter{\openup\jot\m@th
  \ialign{\strut\hfil$\displaystyle{##}$&$\displaystyle{{}##}$\hfil
      \crcr#1\crcr}}\,}
\newif\ifdt@p
\def\displ@y{\global\dt@ptrue\openup\jot\m@th
  \everycr{\noalign{\ifdt@p \global\dt@pfalse \ifdim\prevdepth>-1000\p@
      \vskip-\lineskiplimit \vskip\normallineskiplimit \fi
      \else \penalty\interdisplaylinepenalty \fi}}}
\def\@lign{\tabskip\z@skip\everycr{}} 
\def\displaylines#1{\displ@y \tabskip\z@skip
  \halign{\hbox to\displaywidth{$\@lign\hfil\displaystyle##\hfil$}\crcr
    #1\crcr}}
\def\eqalignno#1{\displ@y \tabskip\plaincentering
  \halign to\displaywidth{\hfil$\@lign\displaystyle{##}$\tabskip\z@skip
    &$\@lign\displaystyle{{}##}$\hfil\tabskip\plaincentering
    &\llap{$\@lign##$}\tabskip\z@skip\crcr
    #1\crcr}}
\def\leqalignno#1{\displ@y \tabskip\plaincentering
  \halign to\displaywidth{\hfil$\@lign\displaystyle{##}$\tabskip\z@skip
    &$\@lign\displaystyle{{}##}$\hfil\tabskip\plaincentering
    &\kern-\displaywidth\rlap{$\@lign##$}\tabskip\displaywidth\crcr
    #1\crcr}}
\def\plaincases#1{\left\{\,\vcenter{\normalbaselines\m@th
    \ialign{$##\hfil$&\quad##\hfil\crcr#1\crcr}}\right.}
\def\plainmatrix#1{\null\,\vcenter{\normalbaselines\m@th
    \ialign{\hfil$##$\hfil&&\quad\hfil$##$\hfil\crcr
      \mathstrut\crcr\noalign{\kern-\baselineskip}
      #1\crcr\mathstrut\crcr\noalign{\kern-\baselineskip}}}\,}

\catcode`\@=12

\def\dlraw{\mathrel{\rlap{$\longrightarrow$}\kern-1pt\longrightarrow}}
\def\vlra{\mathrel{\smash-}\joinrel\mathrel{\smash-}\joinrel%
\kern-2pt\longrightarrow}
\def\srelbar{\vrule width0.6ex height0.65ex depth-0.55ex}
\def\merto{\mathrel{\srelbar\kern1.3pt\srelbar\kern1.3pt\srelbar
    \kern1.3pt\srelbar\kern-1ex\raise0.28ex\hbox{${\scriptscriptstyle>}$}}}

\newdimen\claimskip \claimskip=7pt
\long\def\claim#1|#2\endclaim
{\removelastskip\vskip\claimskip\noindent{\bf#1.}
{\it\ignorespaces#2}\vskip\claimskip\noindent}

\font\ninerm=cmr9
\font\ninei=cmmi9
\font\ninesy=cmsy9
\font\ninebf=cmbx9
\font\ninett=cmtt9
\font\nineit=cmti9
\font\ninesl=cmsl9
\font\ninebb=msbm10 at 9pt
\font\eightrm=cmr8
\font\eighti=cmmi8
\font\eightsy=cmsy8
\font\eightbf=cmbx8

\font\eightit=cmti8
\font\eightsl=cmsl8
\font\sixrm=cmr6
\font\sixi=cmmi6
\font\sixsy=cmsy6
\font\sixbf=cmbx6
\font\fiverm=cmr5
\font\fivei=cmmi5
\font\fivesy=cmsy5
\font\fivebf=cmbx5

\newfam\itfam
\newfam\slfam
\newfam\bffam
\def\eightpoint{\def\rm{\fam0\eightrm}%
\textfont0=\eightrm \scriptfont0=\sixrm \scriptscriptfont0=\fiverm
 \textfont1=\eighti \scriptfont1=\sixi \scriptscriptfont1=\fivei
 \textfont2=\eightsy \scriptfont2=\sixsy \scriptscriptfont2=\fivesy
 \def\it{\fam\itfam\eightit}%
 \textfont\itfam=\eightit
 \def\sl{\fam\slfam\eightsl}%
 \textfont\slfam=\eightsl
 \def\bf{\fam\bffam\eightbf}%
 \textfont\bffam=\eightbf \scriptfont\bffam=\sixbf
 \scriptscriptfont\bffam=\fivebf
 \normalbaselineskip=9pt
 \setbox\strutbox=\hbox{\vrule height7pt depth2pt width0pt}%
 \normalbaselines\rm}

\def\ninepoint{\def\rm{\fam0\ninerm}%
\textfont0=\ninerm \scriptfont0=\sixrm \scriptscriptfont0=\fiverm
 \textfont1=\ninei \scriptfont1=\sixi \scriptscriptfont1=\fivei
 \textfont2=\ninesy \scriptfont2=\sixsy \scriptscriptfont2=\fivesy
 \def\it{\fam\itfam\nineit}%
 \textfont\itfam=\nineit
 \def\sl{\fam\slfam\ninesl}%
 \textfont\slfam=\ninesl
 \def\bf{\fam\bffam\ninebf}%
 \textfont\bffam=\ninebf \scriptfont\bffam=\sixbf
 \scriptscriptfont\bffam=\fivebf
 \normalbaselineskip=11pt
 \setbox\strutbox=\hbox{\vrule height7pt depth2pt width0pt}%
 \normalbaselines\rm}

\def\plainsection#1{\par\vskip .5cm\penalty -100 
\vbox{\noindent{\sc #1}
\vskip 5pt}
\penalty 500}

\def\Bibitem#1&#2&#3&#4&%
{\hangindent=1.66cm\hangafter=1
\noindent\rlap{\hbox{\rm #1}}\kern1.66cm{\rm #2}{\it #3}{\rm #4.}
\vskip5pt}

\def\square{{\hfill \hbox{
\vrule height 1.453ex  width 0.093ex  depth 0ex
\vrule height 1.5ex  width 1.3ex  depth -1.407ex\kern-0.1ex
\vrule height 1.453ex  width 0.093ex  depth 0ex\kern-1.35ex
\vrule height 0.093ex  width 1.3ex  depth 0ex}}}
\def\bigsquare{{\kern-0.3ex\hbox{
\vrule height 1.7ex  width 0.093ex  depth 0ex\kern-0.093ex
\vrule height 1.8ex  width 1.7ex  depth -1.707ex\kern-0.093ex
\vrule height 1.7ex  width 0.093ex  depth 0ex\kern-1.65ex
\vrule height 0.093ex  width 1.6ex  depth 0ex}\kern0.3ex}}
\def\qed{\phantom{~}$\square$\medskip}
\def\smallskip{\vskip 3pt}
\def\medskip{\vskip 5pt}

\title[On the hyperbolicity of very general hypersurfaces]{Proof of the Kobayashi conjecture on the hyperbolicity of very general hypersurfaces}

\author{Jean-Pierre Demailly}
\date{January 23, 2015, revised on March 11, 2015}
\begin{document}

\begin{abstract}
The Green-Griffiths-Lang conjecture stipulates that for every projective 
variety $X$ of general type over $\bC$, there exists a proper algebraic 
subvariety of $X$ containing all non constant entire curves $f:\bC\to X$. 
Using the formalism of directed varieties, we prove here that this 
assertion holds true in case $X$ satisfies a strong general type condition
that is related to a certain jet-semistability property of the tangent
bundle~$T_X$. We then use this fact to confirm a long-standing conjecture
of Kobayashi (1970), according to which a very general algebraic hypersurface
of dimension $n$ and degree at least $2n+2$ in the complex projective space
$\bP^{n+1}$ is hyperbolic.
\end{abstract}

\maketitle
\vskip10pt
\hbox to \textwidth{\hfill\it dedicated to the memory of Salah Baouendi}
\vskip20pt

\plainsection{0. Introduction} 

The goal of this paper, among other results, is to prove the long
standing conjecture of Kobayashi [Kob70, Kob78], according to which a 
very general algebraic hypersurface of dimension $n$ and degree $d\ge 2n+2$ 
in complex projective space $\bP^{n+1}$ is Kobayashi hyperbolic. It is
expected that the bound can be improved to $2n+1$ for $n\ge 2$, and
such a bound would be optimal by Zaidenberg [Zai87], but we cannot yet
prove this. Siu [Siu02, Siu04, Siu12] has introduced a more explicit but 
more computationally involved approach that yields the same conclusion for
$d\ge d_n$, with a very large bound $d_n$ instead of~$2n+2$. However,
thanks to famous results of Clemens [Cle86], Ein [Ein88, Ein91] and
Voisin [Voi96, Voi98], it was known that the bound $2n+2$ would be a
consequence of the Green-Griffiths-Lang conjecture on entire curve
loci, cf.~[GG79] and [Lan86]. Our technique consists in studying a 
generalized form of the GGL conjecture, and proving a special case
that is strong enough to imply the Kobayashi conjecture, using e.g.\
[Voi96]. For this purpose, as was already observed in [Dem97], it is
useful to work in the category of directed projective varieties, and
to take into account the singularities that may appear in the directed
structures, at all steps of the proof.

Since the basic problems we deal with are birationally invariant, the
varieties under consideration can always be replaced by nonsingular
models. A directed projective manifold is a pair $(X,V)$ where $X$ is
a projective manifold equipped with an analytic linear subspace
$V\subset T_X$, i.e.\ a closed irreducible complex analytic subset $V$
of the total space of~$T_X$, such that each fiber $V_x=V\cap T_{X,x}$
is a complex vector space [If $X$ is not irreducible, $V$ should
rather be assumed to be irreducible merely over each component of $X$,
but we will hereafter assume that our varieties are irreducible]. A
morphism $\Phi:(X,V)\to(Y,W)$ in the category of directed manifolds is
an analytic map $\Phi:X\to Y$ such that $\Phi_*V\subset W$. We refer
to the case $V=T_X$ as being the {\it absolute case}, and to the case
$V=T_{X/S}=\Ker d\pi$ for a fibration $\pi:X\to S$, as being the {\it
  relative case}; $V$ may also be taken to be the tangent space to the
leaves of a singular analytic foliation on~$X$, or maybe even a non
integrable linear subspace of $T_X$.

We are especially interested in {\it entire curves} that are tangent to $V$,
namely non constant holomorphic morphisms $f:(\bC,T_\bC)\to (X,V)$ of directed manifolds. In the absolute case, these are just arbitrary entire curves $f:\bC\to X$. The Green-Griffiths-Lang conjecture, in its strong form, stipulates 

\claim 0.1. GGL conjecture|Let $X$ be a projective variety of general type. Then there exists a proper algebraic variety $Y\subsetneq X$ such that every entire curve $f:\bC\to X$ satisfies $f(\bC)\subset Y$.
\endclaim

\noindent [The weaker form would state that entire curves are
algebraically degenerate, so that\break
$f(\bC)\subset Y_f\subsetneq X$ where $Y_f$ might depend on $f\,$]. 
The smallest admissible algebraic set $Y\subset X$ is by definition the 
{\it entire curve locus} of $X$, defined as the Zariski closure
$$
\ECL(X)=\overline{\bigcup_f f(\bC)}^{\rm Zar}.\leqno(0.2)
$$
If $X\subset\bP^N_\bC$ is defined over a number field $\bK_0$ (i.e.\ by
polynomial equations with equations with coefficients in $\bK_0$) and
$Y=\ECL(X)$, it is expected that for every number field $\bK\supset\bK_0$
the set of $\bK$-points in $X(\bK)\smallsetminus Y$ is
finite, and that this property characterizes $\ECL(X)$ as the smallest
algebraic subset $Y$ of $X$ that has the above property for 
all $\bK$ ([Lan86]). This conjectural arithmetical statement would
be a vast generalization of the Mordell-Faltings theorem, and is 
one of the strong motivations to study the geometric GGL conjecture
as a first step.

\claim 0.3. Problem (generalized GGL conjecture)|Let $(X,V)$ be a 
projective directed manifold. Find
geometric conditions on $V$ ensuring that all entire curves $f:(\bC,
T_\bC)\to(X,V)$ are contained in a proper algebraic subvariety
$Y\subsetneq X$. Does this hold when $(X,V)$ is of general type,
in the sense that the canonical sheaf $K_V$ is big~$?$
\endclaim

As above, we define the entire curve locus set of a pair $(X,V)$ to be 
the smallest admissible algebraic set $Y\subset X$ in the above problem, i.e.
$$
\ECL(X,V)=\overline{\mathop{\kern-40pt\bigcup}_{f:(\bC,T_\bC)\to(X,V)}f(\bC)}^{\rm Zar}.\leqno(0.4)
$$
We say that $(X,V)$ is {\it Brody hyperbolic} if $\ECL(X,V)=\emptyset\,$; as is 
well-known, this is equivalent to Kobayashi hyperbolicity whenever
$X$ is compact.

In case $V$ has no singularities, the {\it canonical sheaf} $K_V$ is 
defined to be $(\det \cO(V))^*$ where $\cO(V)$ is the sheaf of holomorphic
sections of $V$, but in general this naive definition would not work.
Take for instance a generic pencil of
elliptic curves $\lambda P(z)+\mu Q(z)=0$ of degree $3$ in $\bP_\bC^2$,
and the linear space $V$ consisting of the tangents to the fibers
of the rational map $\bP_\bC^2\merto\bP_\bC^1$
defined by $z\mapsto Q(z)/P(z)$. Then $V$ is given by
$$
0\longrightarrow \cO(V)\longrightarrow\cO(T_{\bP_\bC^2})
\build\llra{8ex}^{PdQ-QdP}_{}\cO_{\bP^2_\bC}(6)\otimes\cJ_S
\longrightarrow 0
$$
where $S=\Sing(V)$ consists of the 9 points
$\{P(z)=0\}\cap\{Q(z)=0\}$, and $\cJ_S$ is the corresponding ideal
sheaf of~$S$. Since $\det\cO(T_{\bP^2})=\cO(3)$, we see that
$(\det(\cO(V))^*=\cO(3)$ is ample, thus Problem 0.3 would not have a
positive answer (all leaves are elliptic or singular rational curves
and thus covered by entire curves). An even more ``degenerate''
example is obtained with a generic pencil of conics, in which case
$(\det(\cO(V))^*=\cO(1)$ and $\# S=4$.

If we want to get a positive answer to Problem 0.3, it is
therefore indispensable to give a definition of $K_V$ that incorporates
in a suitable way the singularities of $V\,;$ this will be done in
Def.~1.1 (see also Prop.~1.2). The goal is then to give a positive answer 
to Problem~0.3 under some possibly more restrictive conditions for the 
pair $(X,V)$. These conditions will be expressed in terms of the tower
of Semple jet bundles
$$
(X_k,V_k)\to(X_{k-1},V_{k-1})\to\ldots\to(X_1,V_1)\to(X_0,V_0):=(X,V)
\leqno(0.5)
$$
which we define more precisely in Section~1, following [Dem95]. 
It is constructed inductively by setting $X_k=P(V_{k-1})$
(projective bundle of {\it lines} of $V_{k-1}$), and all $V_k$ have the 
same rank $r=\rank V$, so that $\dim X_k=n+k(r-1)$ where $n=\dim X$.
If $\cO_{X_k}(1)$ is the tautological line bundle over $X_k$ associated
with the projective structure and $\pi_{k,\ell}:X_k\to X_\ell$ 
is the natural projection from $X_k$ to $X_\ell$, $0\le\ell\le k$,
we define the $k$-stage Green-Griffiths locus of $(X,V)$ to be
$$
\GG_k(X,V)=\overline{(X_k\smallsetminus\Delta_k)\cap
\bigcap_{m\in\bN}\left(\hbox{base locus of }\cO_{X_k}(m)\otimes 
\pi_{k,0}^*A^{-1}\right)}
\leqno(0.6)
$$
where $A$ is any ample line bundle on $X$ and $\Delta_k=\bigcup_{2\le \ell\le k}
\pi_{k,\ell}^{-1}(D_\ell$) is the union of ``vertical divisors'' 
(see section~1; 
the vertical divisors play no role and have to be removed in this context). 
Clearly, $\GG_k(X,V)$ does not depend on the choice of~$A$.
The basic vanishing theorem for entire 
curves (cf.\ [GG79], [SY96] and [Dem95]) asserts that for every entire 
curve $f:(\bC,T_\bC)\to(X,V)$, then its $k$-jet $f_{[k]}:(\bC,T_\bC)\to
(X_k,V_k)$ satisfies
$$
f_{[k]}(\bC)\subset \GG_k(X,V),\quad\hbox{hence}\quad
f(\bC)\subset \pi_{k,0}\left(\GG_k(X,V)\right).
\leqno(0.7)
$$
(For this, one uses the fact that $f_{[k]}(\bC)$ is not contained in any
component of $\Delta_k$, cf.~[Dem95]). It is therefore natural to define
the global Green-Griffiths locus of $(X,V)$ to be
$$
\GG(X,V)=\smash{\bigcap_{k\in\bN}}\pi_{k,0}\left(\GG_k(X,V)\right).
\leqno(0.8)
$$
By (0.7) we infer that
$$
\ECL(X,V)\subset\GG(X,V).
\leqno(0.9)
$$
The main result of [Dem11] (Theorem~2.37 and Cor.~3.4) implies the 
following useful information:

\claim 0.10. Theorem|Assume that $(X,V)$ is of ``general type'', i.e.\ that 
the canonical sheaf $K_V$ is big on $X$.
Then there exists an integer $k_0$ such that $\GG_k(X,V)$ is a proper algebraic
subset of $X_k$ for $k\ge k_0$ $[\,$though $\pi_{k,0}(\GG_k(X,V))$ might still 
be equal to $X$ for all $k\,]$.
\endclaim

In fact, if $F$ is an invertible sheaf on $X$ such that
$K_V\otimes F$ is big, the probabilistic estimates of [Dem11, Cor.~2.38 and 
Cor.~3.4] produce sections of
$$
\cO_{X_k}(m)\otimes\pi_{k,0}^*\cO\Big(\frac{m}{kr}\Big(
1+\frac{1}{2}+\ldots+\frac{1}{k}\Big)F\Big)\leqno(0.11)
$$
for $m\gg k\gg 1$. The (long and involved) proof uses a curvature computation
and singular holomorphic Morse inequalities to show that the line bundles 
involved in (0.11) are big on $X_k$ for $k\gg 1$. One applies this to 
$F=A^{-1}$ with $A$ ample on $X$ to produce sections and conclude that 
$\GG_k(X,V)\subsetneq X_k$. 

Thanks to (0.9), the GGL conjecture is satisfied whenever 
$\GG(X,V)\subsetneq X$. By [DMR10], this happens for instance in the
absolute case when $X$ is a generic hypersurface of 
degree~\hbox{$d\ge\smash{2^{n^5}}$ in $\bP^{n+1}$} (see also [Pau08], e.g.\
for better bounds in low dimensions).
However, as already mentioned in [Lan86], very simple examples show that 
one can have $\GG(X,V)=X$ even when $(X,V)$ is of general type, and 
this already occurs in the absolute case as soon as $\dim X\ge 2$. 
A typical example is a product of directed manifolds
$$(X,V)=(X',V')\times(X'',V''),\qquad
V=\pr^{\prime\,*}V'\oplus\pr^{\prime\prime\,*}V''.\leqno(0.12)$$
The absolute case $V=T_X$, $V'=T_{X'}$, $V''=T_{X''}$ on a product of curves
is the simplest instance. It is then easy to check that $\GG(X,V)=X$, cf.\
(3.2). Diverio and Rousseau [DR13] have given many more such examples, 
including the case of indecomposable varieties $(X,T_X)$, e.g.\ Hilbert modular 
surfaces, or more generally compact quotients of bounded symmetric domains
of rank${}\ge 2$. The problem here is the failure of some sort of stability 
condition that is introduced in Section 3. This leads to a somewhat
technical concept of more manageable directed pairs $(X,V)$ that we call 
{\it strongly of general type}, see Def.~3.1. Our main result can be stated

\claim 0.13. Theorem (partial solution to the generalized GGL conjecture)|Let 
$(X,V)$ be a directed pair that is strongly of general type.
Then the Green-Griffiths-Lang conjecture holds true for~$(X,V)$, namely
$\ECL(X,V)$ is a proper algebraic subvariety of $X$.
\endclaim

The proof proceeds through a complicated induction on $n=\dim X$ and 
$k=\rank V$, which is the main reason why we have to introduce
directed varieties, even in the absolute case. An interesting feature of 
this result is that the conclusion
on $\ECL(X,V)$ is reached without having to know anything about 
the Green-Griffiths locus $\GG(X,V)$, even a posteriori. Nevetherless,
this is not yet enough to confirm the GGL conjecture. Our hope is
that pairs $(X,V)$ that are of general type without being strongly
of general type -- and thus exhibit some sort of ``jet-instability'' -- can
be investigated by different methods, e.g.\ by the diophantine
approximation techniques of McQuillan [McQ98]. However, Theorem~0.13
is strong enough to imply the Kobayashi conjecture on generic hyperbolicity,
thanks to the following concept of algebraic jet-hyperbolicity.

\claim 0.14. Definition|A directed variety $(X,V)$ will be said to be 
algebraically jet-hyperbolic if the induced directed variety structure
$(Z,W)$ on every irreducible algebraic variety $Z$ of~$X$ such that
$\rank W\ge 1$ has a desingularization that is strongly of general 
type $[$see Sections~$2$ and $4$ for the definition of induced directed 
structures and further details$]$. We also say that a projective
manifold $X$ is algebraically jet-hyperbolic if $(X,T_X)$ is.
\endclaim

In this context, Theorem 0.13 yields the following connection between
algebraic jet-hyperbolicity and the analytic concept of Kobayashi 
hyperbolicity.

\claim 0.15. Theorem|Let $(X,V)$ be a directed variety structure on a
projective manifold $X$. Assume that $(X,V)$ is algebraically jet-hyperbolic.
Then $(X,V)$ is Kobayashi hyperbolic.
\endclaim

This strong link is useful to deal with generic hyperbolicity, i.e.\ the 
hyperbolicity of very general fibers in a deformation $\pi:\cX\to S$
(by ``very general fiber'', we mean here a fiber $X_t=\pi^{-1}(t)$ where
$t$ is taken in a complement $S\smallsetminus\bigcup S_\nu$ of a countable
union of algebraic subsets $S_\nu\subsetneq S$). In Section~5, we apply
the above results to analyze the ``relative'' Semple tower of $(\cX,T_{\cX/S})$
versus the ``absolute'' one derived from $(\cX,T_\cX)$, and obtain
in this way:

\claim 0.16. Theorem|Let $\pi:\cX\to S$ be a deformation of complex
projective nonsingular varieties $X_t=\pi^{-1}(t)$ over a smooth irreducible
quasi-projective base $S$. Let $n=\dim X_t$ be the relative dimension
and let $N=\dim S$. Assume that for all $q=N+1,\ldots, N+n$, the
exterior power $\Lambda^qT^*_{\cX}$ is a relatively ample vector
bundle over~$S$. Then the very general fiber $X_t$ is algebraically
jet-hyperbolic, and thus Kobayashi hyperbolic.
\endclaim

In the special case of the universal family of complete intersections
of codimension $c$ and type $(d_1,\ldots,d_c)$ in complex projective
$\bP^{n+c}$, Ein [Ein88, 91] and Voisin [Voi96] have exploited in a crucial 
way the existence of global twisted vector fields, e.g.\ sections
of $T_{\cX}\otimes\cO(1)$, on the total space~$\cX$ over $S$. This idea
combined with Theorem~0.16 yields the following result.

\claim 0.17. Corollary  (confirmation of the Kobayashi conjecture)|If 
$\sum d_j\ge 2n+c+1$, the very general complete
intersection of type $(d_1,\ldots,d_c)$ in complex projective space
$\bP^{n+c}$ is Kobayashi hyperbolic.
\endclaim

The border case $\sum d_j=2n+c$ is potentially also accessible by the 
techniques developed here, but further calculations would be needed to
check the possible degenerations of the morphisms induced by twisted
vector fields. 
I would like to thank Simone Diverio and Erwan Rousseau for very
stimulating discussions on these questions. I am grateful to Mihai
P\u{a}un for an invitation at KIAS (Seoul) in August 2014, during
which further very fruitful exchanges took place, and for his
extremely careful reading of earlier drafts of the manuscript.

\plainsection{1. Semple jet bundles and associated canonical sheaves} 

Let $(X,V)$ be a directed projective manifold and $r=\rank V$, that is, the
dimension of generic fibers.
Then $V$ is actually a holomorphic subbundle of $T_X$ on the complement
$X\smallsetminus \Sing(V)$ of a certain minimal analytic set
$\Sing(V)\subsetneq X$ of codimension${}\ge 2$, called hereafter the 
singular set of~$V$. If $\mu:\smash{\widehat X}\to X$ is a proper modification
(a composition of blow-ups with smooth centers, say), we get a directed
manifold $\smash{(\widehat X,\widehat V)}$ by taking $\smash{\widehat V}$ 
to be the closure of $\mu_*^{-1}(V')$, where $V'=V_{|X'}$ is the
restriction of $V$ over a Zariski open set $X'\subset X\smallsetminus \Sing(V)$ 
such that $\mu:\mu^{-1}(X')\to X'$ is a biholomorphism. We will be 
interested in taking modifications realized by iterated blow-ups of 
certain nonsingular subvarieties of the singular set $\Sing(V)$, so 
as to eventually ``improve'' the singularities of $V\,$; outside of
$\Sing(V)$ the effect of blowing-up will be irrelevant, as one
can see easily. Following [Dem11], the canonical sheaf $K_V$ is defined
as follows.

\claim 1.1. Definition|For any directed pair $(X,V)$ with $X$ nonsingular,
we define $K_V$ to be the rank~$1$ analytic sheaf such that
$$
K_V(U)=\hbox{sheaf of locally bounded sections 
of}~~\cO_X(\Lambda^r V^{\prime *})(U\cap X')
$$
where $r=\hbox{rank}(V)$, $X'=X\smallsetminus \Sing(V)$, $V'=V_{|X'}$,
and ``bounded'' means bounded with respect to a smooth hermitian metric
$h$ on $T_X$.
\endclaim

For $r=0$, one can set $K_V=\cO_X$, but this case is trivial: clearly
$\ECL(X,V)=\emptyset$.  The above definition of $K_V$ may look like an
analytic one, but it can easily be turned into an equivalent algebraic
definition:

\claim 1.2. Proposition|Consider the natural morphism 
$\cO(\Lambda^rT_X^*)\to \cO(\Lambda^r V^*)$ where $r=\rank V$
$[\cO(\Lambda^r V^*)$ being defined here as the quotient of 
$\cO(\Lambda^rT_X^*)$ by $r$-forms that have zero restrictions 
to $\cO(\Lambda^rV^*)$ on $X\smallsetminus \Sing(V)\,]$. The 
bidual $\cL_V=\cO_X(\Lambda^r V^*)^{**}$ is an invertible sheaf,
and our natural morphism can be written
$$
\cO(\Lambda^rT_X^*)\to \cO(\Lambda^rV^*)=\cL_V\otimes\cJ_V\subset \cL_V
\leqno(1.2.1)
$$
where $\cJ_V$ is a certain ideal sheaf of $\cO_X$ whose zero set is
contained in $\Sing(V)$ and the arrow on the left is surjective by 
definition. Then 
$$
K_V=\cL_V\otimes \overline{\cJ}_V\leqno(1.2.2)
$$
where $\overline{\cJ}_V$ is the integral closure of $\cJ_V$ in $\cO_X$.
In particular, $K_V$ is always a coherent sheaf.
\endclaim

\proof Let $(u_k)$ be a set of generators of $\cO(\Lambda^rV^*)$ obtained
(say) as the images of a basis $(dz_I)_{|I|=r}$ of $\Lambda^rT^*_X$
in some local coordinates near a point $x\in X$. Write $u_k=g_k\ell$ 
where $\ell$ is a local generator of $\cL_V$ at~$x$. Then $\cJ_V=(g_k)$ by
definition. The boundedness condition expressed in Def.~1.1 means that we 
take sections of the form $f\ell$ where $f$ is a holomorphic function 
on~$U\cap X'$ (and $U$ a neighborhood of $x$), such that 
$$
|f|\le C\sum|g_k|\leqno(1.2.3)
$$
for some constant $C>0$. But then $f$ extends holomorphically to $U$ into 
a function that lies in the integral closure $\overline\cJ_V$, and the latter
is actually characterized analytically by condition~(1.2.3). This proves
Prop.~1.2.\qed

By blowing-up $\cJ_V$ and taking a desingularization $\widehat X$, one 
can always find a {\it log-resolution} of $\cJ_V$ (or $K_V$), i.e.\ a 
modification $\mu:\smash{\widehat X}\to X$ such 
that $\mu^*\cJ_V\subset\cO_{\widehat X}$ is an invertible ideal
sheaf (hence integrally closed); it~follows that
$\mu^*\overline{\cJ}_V=\mu^*\cJ_V$ and
$\mu^*K_V=\mu^*\cL_V\otimes\mu^*\cJ_V$ are invertible sheaves
on~$\smash{\widehat X}$. Notice that for any modification 
$\mu':(X',V')\to (X,V)$, there is always a well defined
natural morphism
$$
\mu^{\prime\,*}K_V\to K_{V'}\leqno(1.3)
$$
(though it need not be an isomorphism, and $K_{V'}$ is possibly non
invertible even when $\mu'$ is taken to be a log-resolution of
$K_V$). Indeed $(\mu')_*=d\mu':V'\to \mu^*V$ is continuous
with respect to ambient hermitian metrics on $X$ and $X'$, and going
to the duals reverses the arrows while preserving boundedness with
respect to the metrics. If $\mu'':X''\to X'$ provides a simultaneous
log-resolution of $K_{V'}$ and $\mu^{\prime\,*}K_V$, we get a non trivial
morphism of invertible sheaves
$$
(\mu'\circ\mu'')^*K_V= \mu^{\prime\prime\,*}\mu^{\prime\,*}K_V
\longrightarrow \mu^{\prime\prime\,*}K_{V'},\leqno(1.4)
$$
hence the bigness of $\mu^{\prime\,*}K_V$ with imply that of 
$\mu^{\prime\prime\,*}K_{V'}$. This is a general principle that we would like
to refer to as the ``monotonicity principle'' for canonical sheaves:
one always get more sections by going to a higher level through a
(holomorphic) modification.

\claim 1.5. Definition|We say that the rank $1$ sheaf
$K_V$ is ``big'' if the invertible sheaf $\mu^*K_V$ is big 
in the usual sense for any log resolution 
$\mu:\smash{\widehat X}\to X$ of~$K_V$. 
Finally, we say that $(X,V)$ is of {\it general type} if there exists
a modification $\mu':(X',V')\to (X,V)$
such that $K_{\smash{V'}}$ is big$\;;$ any higher blow-up 
$\mu'':(X'',V'')\to (X',V')$ then also yields a big canonical sheaf
by $(1.3)$.
\endclaim

Clearly, ``general type'' is a birationally (or bimeromorphically)
invariant concept, by the very definition.  When
$\dim X=n$ and $V\subset T_X$ is a subbundle of rank $r\ge 1$,
one constructs a tower of ``Semple $k$-jet bundles''
$\pi_{k,k-1}:(X_k,V_k)\to (X_{k-1},V_{k-1})$ that are 
$\bP^{r-1}$-bundles, with \hbox{$\dim X_k=n+k(r-1)$} and
$\rank(V_k)=r$. For this, we take $(X_0,V_0)=(X,V)$, and for every 
$k\ge 1$, we set inductively $X_k:=P(V_{k-1})$ and
$$
V_k:=(\pi_{k,k-1})_*^{-1}\cO_{X_k}(-1)\subset T_{X_k},
$$
where $\cO_{X_k}(1)$ is the tautological line bundle on $X_k$,
$\pi_{k,k-1}:X_k=P(V_{k-1})\to X_{k-1}$ the natural projection and
$(\pi_{k,k-1})_*=d\pi_{k,k-1}:T_{X_k}\to
\pi^*_{k,k-1}T_{X_{k-1}}$ its differential (cf.\ [Dem95]). In other terms,
we have exact sequences
$$
\leqalignno{
&0\longrightarrow T_{X_k/X_{k-1}}\longrightarrow V_k
\mathop{\longrightarrow}\limits^{(\pi_{k,k-1})_*}\cO_{X_k}(-1)
\longrightarrow 0,&(1.6)\cr
&0\longrightarrow \cO_{X_k}\longrightarrow(\pi_{k,k-1})^*V_{k-1}\otimes
\cO_{X_k}(1)
\longrightarrow T_{X_k/X_{k-1}}\longrightarrow 0,&(1.7)\cr}
$$
where the last line is the Euler exact sequence associated with the 
relative tangent bundle of $P(V_{k-1})\to X_{k-1}$. Notice that we
by definition of the tautological line bundle we have
$$
\cO_{X_k}(-1)\subset\pi_{k,k-1}^*V_{k-1}
\subset\pi_{k,k-1}^*T_{X_{k-1}},
$$
and also $\rank(V_k)=r$. Let us recall also that for $k\ge 2$, there 
are ``vertical divisors''
$D_k=P(T_{X_{k-1}/X_{k-2}})\subset P(V_{k-1})=X_k$, and that
$D_k$ is the zero divisor of the section of $\cO_{X_k}(1)\otimes
\pi_{k,k-1}^*\cO_{X_{k-1}}(-1)$ induced by the second arrow of the
first exact sequence (1.6), when $k$ is replaced by $k-1$. 
This yields in particular
$$
\cO_{X_k}(1)=\pi_{k,k-1}^*\cO_{X_{k-1}}(1)\otimes\cO(D_k).\leqno(1.8)
$$
By composing the projections we get for all pairs of indices 
$0\le j\le k$ natural morphisms
$$
\pi_{k,j}:X_k\to X_j,\quad
(\pi_{k,j})_*=(d\pi_{k,j})_{|V_k}:V_k\to(\pi_{k,j})^*V_j,
$$
and for every $k$-tuple ${\bf a}=(a_1,\ldots,a_k)\in\bZ^k$
we define
$$
\cO_{X_k}({\bf a})=\bigotimes_{1\le j\le k}\pi_{k,j}^*\cO_{X_j}(a_j),\quad
\pi_{k,j}:X_k\to X_j.
$$
We extend this definition to all weights ${\bf a}\in\bQ^k$ to get
a $\bQ$-line bundle in $\Pic(X)\otimes_\bZ\bQ$. Now, Formula (1.8) yields
$$
\cO_{X_k}({\bf a})=\cO_{X_k}(m)\otimes \cO(-{\bf b}\cdot D)\quad
\hbox{where~~$m=|{\bf a}|=\sum a_j$, ${\bf b}=(0,b_2,\ldots,b_k)$}
\leqno(1.9)
$$
and $b_j=a_1+\ldots+a_{j-1}$, $2\le j\le k$.

When $\Sing(V)\neq\emptyset$, one can always
define $X_k$ and $V_k$ to be the respective closures of $X'_k$, $V'_k$ 
associated with $X'=X\smallsetminus \Sing(V)$ and $V'=V_{|X'}$, where
the closure is taken in the nonsingular ``absolute''
Semple tower $(X^a_k,V^a_k)$ obtained from $(X^a_0,V^a_0)=(X,T_X)$. 
We leave the reader check the following easy (but important) observation.

\claim 1.10. Fonctoriality|If $\Phi:(X,V)\to(Y,W)$ is a morphism of
directed varieties such that $\Phi_*:T_X\to\Phi^*T_Y$ is injective
$($i.e.\ $\Phi$ is an immersion$\,)$, then there is a corresponding
natural morphism $\Phi_{[k]}:(X_k,V_k)\to(Y_k,W_k)$
at the level of Semple bundles. If one merely assumes that the differential
$\Phi_*:V\to\Phi^*W$ is non zero, there is still a well defined
meromorphic map $\Phi_{[k]}:(X_k,V_k)\merto(Y_k,W_k)$ for all~$k\ge 0$.
\endclaim

In case $V$ is singular, the $k$-th Semple bundle $X_k$ will also
be singular, but we can still replace $(X_k,V_k)$ by a suitable modification 
$(\smash{\widehat X}_k,\smash{\widehat V}_k)$ if we want to work with 
a nonsingular model $\smash{\widehat X}_k$ of~$X_k$. The exceptional
set of $\smash{\widehat X}_k$ over $X_k$ can be chosen to lie above
$\Sing(V)\subset X$, and proceeding inductively with respect to $k$, we can 
also arrange the modifications in such a way that we get a tower 
structure $\smash{({\widehat X}_{k+1},
{\widehat V}_{k+1})}\to\smash{({\widehat X}_k,{\widehat V}_k)}\,$; however,
in general, it will not be possible to achieve that 
$\smash{{\widehat V}_k}$ is a subbundle of $T_{\widehat X_k}$.

It is not true that
$K_{\smash{\widehat V}_k}$ is big in case $(X,V)$ is of general type (especially
since the fibers of $X_k\to X$ are towers of $\bP^{r-1}$ bundles, and the
canonical bundles of projective spaces are always negative~!). However,
a twisted version holds true, that can be seen as another instance
of the ``monotonicity principle'' when going to higher stages
in the Semple tower.

\claim 1.11. Lemma|If $(X,V)$ is of general type, then there is a
modification $\smash{(\widehat X,\widehat V)}$ such that all pairs
$\smash{(\widehat X_k,\widehat V_k)}$ of the associated Semple tower
have a twisted canonical bundle $\smash{K_{\widehat V_k}\otimes
\cO_{\widehat X_k}(p)}$ that is still big when one multiplies $K_{\widehat V_k}$ by 
a suitable $\bQ$-line bundle $\cO_{\widehat X_k}(p)$, $p\in\bQ_+$.
\endclaim

\proof First assume that $V$ has no singularities.
The exact sequences (1.6) and (1.7) provide
$$
K_{V_k}:=\det V_k^*=\det(T^*_{X_k/X_{k-1}})\otimes\cO_{X_k}(1)=
\pi^*_{k,k-1}K_{V_{k-1}}\otimes\cO_{X_k}(-(r-1))
$$
where $r=\hbox{rank}(V)$. Inductively we get
$$
K_{V_k}=\pi_{k,0}^*K_V\otimes\cO_{X_k}(-(r-1){\bf 1}),\qquad
{\bf 1}=(1,...,1)\in\bN^k.\leqno(1.11.1)
$$
We know by [Dem95] that $\cO_{X_k}({\bf c})$ is relatively ample over
$X$ when we take the special weight
${\bf c}=(2\,3^{k-2},...,2\,3^{k-j-1},...,6,2,1)$, hence
$$
K_{V_k}\otimes\cO_{X_k}((r-1){\bf 1}+\varepsilon{\bf c})
=\pi_{k,0}^* K_V\otimes\cO_{X_k}(\varepsilon{\bf c})
$$
is big over $X_k$ for any sufficiently small positive rational number 
$\varepsilon\in\bQ^*_+$. Thanks to Formula~(1.9), we can in fact
replace the weight $(r-1){\bf 1}+\varepsilon{\bf c}$ by its total degree
$p=(r-1)k+\varepsilon|{\bf c}|\in\bQ_+$. The general case of a singular 
linear space follows by considering suitable ``sufficiently high'' 
modifications $\smash{\widehat X}$ of $X$, the related directed structure
$\smash{\widehat V}$ on $\smash{\widehat X}$, and embedding 
$\smash{(\widehat X_k,\widehat V_k)}$ in the absolute Semple tower
$\smash{(\widehat X^a_k,\widehat V^a_k)}$ of $\smash{\widehat X}$.
We still have a well defined morphism of rank $1$ sheaves
$$
\pi_{k,0}^*K_{\widehat V}\otimes\cO_{\widehat X_k}(-(r-1){\bf 1})\to
K_{\widehat V_k}
\leqno(1.11.2)
$$
because the multiplier ideal sheaves involved at each stage behave according
to the monoto\-nicity principle applied to the projections
$\smash{\pi^a_{k,k-1}:\widehat X^a_k\to \widehat X^a_{k-1}}$ and their
differentials $(\pi^a_{k,k-1})_*$, which yield well-defined transposed
morphisms from the $(k-1)$-st stage to the $k$-th stage at the level of
exterior differential forms. Our contention follows.\qed

\plainsection{2. Induced directed structure on a subvariety
of a jet space}

Let $Z$ be an irreducible algebraic subset of some $k$-jet bundle
$X_k$ over~$X$, $k\ge 0$. We define the linear subspace 
\hbox{$W\subset T_Z\subset T_{X_k|Z}$} to be the closure
$$
W:=\overline{T_{Z'}\cap V_k}\leqno(2.1)
$$
taken on a suitable Zariski open set $Z'\subset Z_{\rm reg}$ where 
the intersection $T_{Z'}\cap V_k$ has constant rank and 
is a subbundle of $T_{Z'}$. Alternatively, we could also take 
$W$ to be the closure of $T_{Z'}\cap V_k$ in the $k$-th stage
$(X^a_k,V^a_k)$ of the absolute Semple tower. We say that $(Z,W)$ is the 
{\it induced} directed variety structure. In the sequel, we always
consider such a subvariety $Z$ of $X_k$ as a directed pair $(Z,W)$ by
taking the induced structure described above. Let us first quote the 
following easy observation.

\claim 2.2. Observation|For $k\ge 1$, let 
$Z\subsetneq X_k$ be an irreducible algebraic subset that projects
onto~$X_{k-1}$, i.e.\ $\pi_{k,k-1}(Z)=X_{k-1}$.
Then the induced directed variety \hbox{$(Z,W)\subset(X_k,V_k)$},
satisfies
$$1\le \rank W<r:=\rank(V_k).$$
\endclaim

\proof Take a Zariski open subset $Z'\subset Z_{\rm reg}$
such that $W'=T_{Z'}\cap V_k$ is a vector bundle over $Z'$. Since 
$X_k\to X_{k-1}$ is a $\bP^{r-1}$-bundle, $Z$ has codimension at most 
$r-1$ in $X_k$. Therefore $\rank W\ge \rank V_k-(r-1)\ge 1$.
On the other hand, if we had $\rank W=\rank V_k$ generically, then
$T_{Z'}$ would contain $V_{k|Z'}$, in particular it would contain all
vertical directions $T_{X_k/X_{k-1}}\subset V_k$ that are tangent
to the fibers of $X_k\to X_{k-1}$. By taking the flow along vertical 
vector fields, we would
conclude that $Z'$ is a union of fibers of $X_k\to X_{k-1}$ up to an
algebraic set of smaller dimension, but this is excluded since $Z$
projects onto $X_{k-1}$ and $Z\subsetneq X_k$.\qed

\claim 2.3. Definition|For $k\ge 1$, let 
$Z\subset X_k$ be an irreducible algebraic subset of $X_k$ that
projects onto~$X_{k-1}$. We assume moreover that 
$Z\not\subset D_k=P(T_{X_{k-1}/X_{k-2}})$ 
$($and~put here $D_1=\emptyset$ in what follows to avoid to have to 
single out the case $k=1)$. In this situation
we say that $(Z,W)$ is of general type modulo $X_k\to X$ if there exists 
$p\in\bQ_+$ such that $K_W\otimes\cO_{X_k}(p)_{|Z}$ is big over~$Z$, 
possibly after
replacing $Z$ by a suitable nonsingular model $\smash{\widehat Z}$
$($and pulling-back $W$ and $\cO_{X_k}(p)_{|Z}$ to the nonsingular
variety $\smash{\widehat Z}\,)$.
\endclaim

The main result of [Dem11] mentioned in the introduction as 
Theorem~0.10 implies the following important ``induction step''.

\claim 2.4. Proposition|Let $(X,V)$ be a directed pair where
$X$ is projective algebraic. Take an irreducible algebraic subset
$Z\not\subset D_k$ of the associated $k$-jet Semple bundle
$X_k$ that projects onto~$X_{k-1}$, $k\ge 1$, and assume that the 
induced directed space $(Z,W)\subset(X_k,V_k)$ is of general type 
modulo $X_k\to X$.  Then there exists a divisor $\Sigma\subset Z_\ell$
in a sufficiently high stage of the Semple tower $(Z_\ell,W_\ell)$
associated with $(Z,W)$, such that every non constant holomorphic map
$f:\bC\to X$ tangent to $V$ that satisfies $f_{[k]}(\bC)\subset Z$
also satisfies $f_{[k+\ell]}(\bC)\subset\Sigma$.
\endclaim

\proof Let $E\subset Z$ be a divisor containing 
$Z_{\rm sing}\cup(Z\cap \pi_{k,0}^{-1}(\Sing(V)))$, chosen so that
on the nonsingular Zariski open set $Z'=Z\smallsetminus E$ all linear spaces
$T_{Z'}$, $V_{k|Z'}$ and $W'=T_{Z'}\cap V_k$ are subbundles of $T_{X_k|Z'}$,
the first two having a transverse intersection on $Z'$.
By taking closures over $Z'$ in the
absolute Semple tower of $X$, we get (singular) directed pairs 
$(Z_\ell,W_\ell)\subset (X_{k+\ell},V_{k+\ell})$, which we eventually
resolve into 
$(\smash{\widehat Z}_\ell,\smash{\widehat W}_{\ell})\subset 
(\smash{\widehat X}_{k+\ell},\smash{\widehat V}_{k+\ell})$ over
nonsingular bases. By construction, locally bounded sections of 
$\cO_{\smash{\widehat X}_{k+\ell}}(m)$ restrict to locally bounded 
sections of $\cO_{\smash{\widehat Z}_{\ell}}(m)$ over 
$\smash{\widehat Z}_{\ell}$.

Since Theorem~0.10 and the related estimate (0.11) are universal in the 
category of directed varieties, we can apply them by replacing $X$
with $\smash{\widehat Z}\subset \smash{\widehat X}_k$, the order $k$ by
a new index $\ell$, and $F$ by 
$$
F_k=\mu^*\Big(\big(\cO_{X_k}(p)\otimes\pi_{k,0}^*\cO_X(-\varepsilon A)
\big)_{|Z}\Big)
$$
where $\mu:\widehat Z\to Z$ is the desingularization, $p\in\bQ_+$ is
chosen such that $K_W\otimes\cO_{x_k}(p)_{|Z}$ is big, 
$A$ is an ample bundle on $X$ and $\varepsilon\in\bQ_+^*$is  small enough.
The assumptions
show that $K_{\widehat{W}}\otimes F_k$ is big on $\smash{\widehat Z}$,
therefore, by applying our theorem and taking $m\gg\ell\gg 1$, 
we get in fine a large number of (metric bounded) sections of
$$\eqalign{
\cO_{\widehat Z_\ell}(m)&\otimes{\widehat \pi}_{k+\ell,k}^*
\cO\Big(\frac{m}{\ell r'}\Big(
1+\frac{1}{2}+\ldots+\frac{1}{\ell}\Big)F_k\Big)\cr
&=
\cO_{{\widehat X}_{k+\ell}}(m{\bf a'})\otimes{\widehat\pi}_{k+\ell,0}^*\cO\Big(
-\frac{m\varepsilon}{kr}\Big(
1+\frac{1}{2}+\ldots+\frac{1}{k}\Big)A\Big)_{|{\widehat Z}_\ell}\cr}
$$
where ${\bf a'}\in\bQ_+^{k+\ell}$ is a positive weight 
(of the form $(0,\ldots,\lambda,\ldots,0,1)$ with some
non zero component $\lambda\in\bQ_+$ at index $k$).
These sections descend to metric bounded sections of
$$
\cO_{X_{k+\ell}}((1+\lambda)m)\otimes{\widehat\pi}_{k+\ell,0}^*\cO\Big(
-\frac{m\varepsilon}{kr}\Big(
1+\frac{1}{2}+\ldots+\frac{1}{k}\Big)A\Big)_{|Z_\ell}.
$$
Since $A$ is ample on $X$, we can apply the fundamental vanishing theorem 
(see e.g.\ [Dem97] or [Dem11], Statement~8.15), or
rather an ``embedded'' version for curves satis\-fying $f_{[k]}(\bC)\subset Z$,
proved exactly by the same arguments. The vanishing theorem
implies that the divisor $\Sigma$ of any such section satisfies the
conclusions of Proposition~2.4, possibly modulo exceptional divisors
of $\smash{\widehat Z}\to Z$; to take care of these, it is enough to add
to $\Sigma$ the inverse image of the divisor $E=Z\smallsetminus Z'$ 
initially selected.\qed

\plainsection{3. Strong general type condition for directed manifolds}

Our main result is the following partial solution to the Green-Griffiths-Lang
conjecture, providing a sufficient algebraic condition for the analytic
conclusion to hold true. We first give an ad hoc definition.

\claim 3.1. Definition|Let $(X,V)$ be a directed pair where
$X$ is projective algebraic. We say that that $(X,V)$ is ``strongly of
general type'' if it is of general type and for every irreducible 
algebraic set $Z\subsetneq X_k$, $Z\not\subset D_k$, that projects 
onto~$X_{k-1}$, $k\ge 1$, the induced directed structure
$(Z,W)\subset(X_k,V_k)$ is of general type modulo $X_k\to X$.
\endclaim

\claim 3.2. Example|\rm The situation of a product 
$(X,V)=(X',V')\times(X'',V'')$ described in (0.12) shows that $(X,V)$ can
be of general type without being strongly of general type. In fact,
if $(X',V')$ and $(X'',V'')$ are of general type, then 
$K_V=\pr^{\prime\,*}K_{V'}\otimes\pr^{\prime\prime\,*}K_{V''}$ is big, so $(X,V)$ is 
again of general type. However 
$$
Z=P(\pr^{\prime\,*}V')=X'_1\times X''\subset X_1
$$
has a directed structure $W=\pr^{\prime\,*}V'_1$ which does not possess
a big canonical bundle over $Z$, since the restriction of $K_W$ to
any fiber $\{x'\}\times X''$ is trivial. The higher stages $(Z_k,W_k)$
of the Semple tower of $(Z,W)$ are given by $Z_k=X'_{k+1}\times X''$ and
$W_k=\pr^{\prime\,*}V'_{k+1}$, so it is easy to see that $\GG_k(X,V)$ contains
$Z_{k-1}$. Since $Z_k$ projects onto $X$, we have here $\GG(X,V)=X$
(see [DR13] for more sophisticated indecomposable examples).
\endclaim

\claim 3.3. Remark|\rm It follows from Definition~2.3 that 
$(Z,W)\subset(X_k,V_k)$ is automatically of general type modulo $X_k\to X$ 
if $\cO_{X_k}(1)_{|Z}$ is big. Notice further that
$$
\cO_{X_k}(1+\varepsilon)_{|Z}=
\big(\cO_{X_k}(\varepsilon)\otimes\pi_{k,k-1}^*\cO_{X_{k-1}}(1)
\otimes\cO(D_k)\big)_{|Z}
$$
where $\cO(D_k)_{|Z}$ is effective and $\cO_{X_k}(1)$ is relatively
ample with respect to the projection $X_k\to X_{k-1}$. Therefore the
bigness of $\cO_{X_{k-1}}(1)$ on $X_{k-1}$ also implies that every
directed subvariety $(Z,W)\subset(X_k,V_k)$ is of general type modulo
$X_k\to X$. If
$(X,V)$ is of general type, we know by the main result of [Dem11] that
$\cO_{X_k}(1)$
is big for $k\ge
k_0$ large enough, and actually the precise estimates obtained therein
give explicit bounds for such a~$k_0$.
The above observations show that we need to check the condition of
Definition~3.1 only for $Z\subset
X_k$, $k\le k_0$.  Moreover, at least in the case where
$V$,~$Z$, and $W=T_Z\cap V_k$ are nonsingular, we have
$$
K_W\simeq K_Z\otimes\det (T_Z/W)\simeq K_Z\otimes \det(T_{X_k}/V_k)_{|Z}
\simeq K_{Z/X_{k-1}}\otimes\cO_{X_k}(1)_{|Z}.
$$
Thus we see that, in some sense, it is only needed to check the bigness 
of $K_W$ modulo $X_k\to X$ for ``rather special subvarieties'' $Z\subset X_k$ 
over $X_{k-1}$, such that $K_{Z/X_{k-1}}$ is not relatively big over $X_{k-1}$.\qed
\endclaim

\claim 3.4. Hypersurface case|\rm Assume that $Z\ne D_k$ is an irreducible 
hypersurface of $X_k$ that projects onto $X_{k-1}$. To simplify things further,
also assume that $V$ is nonsingular. Since the Semple jet-bundles $X_k$ 
form a tower of $\bP^{r-1}$-bundles, their Picard groups satisfy
$\Pic(X_k)\simeq\Pic(X)\oplus\bZ^k$ and we have 
$\cO_{X_k}(Z)\simeq\cO_{X_k}({\bf a})\otimes\pi_{k,0}^*B$
for some ${\bf a}\in\bZ^k$ and $B\in\Pic(X)$, where $a_k=d>0$ is the
relative degree of the hypersurface over $X_{k-1}$. Let $\sigma\in H^0(X_k,
\cO_{X_k}(Z))$ be the section defining $Z$ in $X_k$.
The induced directed variety $(Z,W)$ has $\rank W=r-1=\rank V-1$ and
formula (1.12) yields $K_{V_k}=\cO_{X_k}(-(r-1){\bf 1})\otimes\pi_{k,0}^*(K_V)$.
We claim that
$$
K_W\supset\big(K_{V_k}\otimes \cO_{X_k}(Z)\big)_{|Z}
\otimes\cJ_S=
\big(\cO_{X_k}({\bf a}-(r-1){\bf 1})\otimes\pi_{k,0}^*(B\otimes K_V)\big)_{|Z}
\otimes\cJ_S\leqno(3.4.1)
$$
where $S\subsetneq Z$ is the set (containing $Z_{\rm sing}$) where $\sigma$
and $d\sigma_{|V_k}$ both vanish, and $\cJ_S$ is the ideal locally 
generated by the coefficients of $d\sigma_{|V_k}$ along $Z=\sigma^{-1}(0)$. In
fact, the intersection $W=T_Z\cap V_k$ is transverse on $Z\smallsetminus S\,$;
then (3.4.1) can be seen by looking at the morphism
$$
V_{k|Z}\build\llra{4ex}^{d\sigma_{|V_k}}_{}\cO_{X_k}(Z)_{|Z},
$$
and observing that the contraction by $K_{V_k}=\Lambda^rV_k^*$ provides a
metric bounded section of the canonical sheaf $K_W$. In order to investigate
the positivity properties of  $K_W$, one has to show that $B$ cannot be too
negative, and in addition to control the singularity set $S$. The second
point is a priori very challenging, but we get useful information for
the first point by observing that $\sigma$ provides a morphism
$\pi_{k,0}^*\cO_X(-B)\to\cO_{X_k}({\bf a})$, hence a nontrivial morphism
$$
\cO_X(-B)\to E_{\bf a}:=(\pi_{k,0})_*\cO_{X_k}({\bf a})
$$
By [Dem95, Section~12], there exists a filtration on $E_{\bf a}$ such that 
the graded pieces are irreducible representations of $\GL(V)$ contained
in $(V^*)^{\otimes \ell}$, $\ell\le|{\bf a}|$. Therefore we get a nontrivial 
morphism
$$
\cO_X(-B)\to (V^*)^{\otimes \ell},\qquad \ell\le|{\bf a}|.\leqno(3.4.2)
$$
If we know about certain (semi-)stability properties of $V$, this can be used
to control the negativity of $B$.\qed
\endclaim

\noindent
We further need the following useful concept that generalizes entire curve loci.

\claim 3.5. Definition|If $Z$ is an algebraic set contained
in some stage $X_k$ of the Semple tower of~$(X,V)$, we define its
``induced entire curve locus'' $\IEL_{X,V}(Z)\subset Z$ to be the Zariski closure 
of the union $\bigcup f_{[k]}(\bC)$ of all jets of entire curves
$f:(\bC,T_\bC)\to(X,V)$ such that $f_{[k]}(\bC)\subset Z$.
\endclaim

We~have of course $\IEL_{X,V}(\IEL_{X,V}(Z))=\IEL_{X,V}(Z)$ by definition.
It is not hard to check that modulo certain ``vertical divisors'' of $X_k$, 
the $\IEL_{X,V}(Z)$ locus is essentially the same as 
the entire curve locus $\ECL(Z,W)$ of the induced directed variety,
but we will not use this fact here. Since $\IEL_{X,V}(X)=\ECL(X,V)$,
proving the Green-Griffiths-Lang property amounts to showing that
$\IEL_{X,V}(X)\subsetneq X$ in the stage $k=0$ of the tower. 

\claim 3.6. Theorem|Let $(X,V)$ be a directed pair of general type. 
Assume that there is an integer $k_0\ge 0$ such that for every $k>k_0$
and every irreducible algebraic set $Z\subsetneq X_k$, $Z\not\subset D_k$, 
that projects onto~$X_{k-1}$, the induced directed structure
$(Z,W)\subset(X_k,V_k)$ is of general type modulo $X_k\to X$.
Then $\IEL_{X,V}(X_{k_0})\subsetneq X_{k_0}$.
\endclaim

\proof We argue here by contradiction,
assuming that $\IEL_{X,V}(X_{k_0})=X_{k_0}$. The main argument consists of producing 
inductively an increasing sequence of integers
$$k_0<k_1<\ldots<k_j<\ldots$$
and directed varieties $(Z^j,W^j)\subset(X_{k_j},V_{k_j})$ satisfying
the following properties~:

{\plainitemindent=14mm
\plainitem{(3.6.1)} $(Z^0,W^0)=(X_{k_0},V_{k_0})\;$;

\plainitem{(3.6.2)} for all $j\ge 0$, $\IEL_{X,V}(Z^j)=Z^j\;$;

\plainitem{(3.6.3)} $Z^j$ is an irreducible algebraic variety
such that $Z^j\subsetneq X_{k_j}$ for $j\ge 1$, 
$Z^j$ is not contained in the vertical divisor $D_{k_j}=
\smash{P(T_{X_{k_j-1}/X_{k_j-2}})}$ of $X_{k_j}$, and
$(Z^j,W^j)$ is of general type modulo $X_{k_j}\to X$
(i.e.\ some nonsingular model is)$\;$;

\plainitem{(3.6.4)} for all $j\ge 0$, the directed variety 
$(Z^{j+1},W^{j+1})$ is contained 
in some stage (of order $\ell_j=k_{j+1}-k_j$) of the Semple
tower of $(Z^j,W^j)$, namely 
$$
(Z^{j+1},W^{j+1})\subset (Z^j_{\ell_j},W^j_{\ell_j})\subset
(X_{k_{j+1}},V_{k_{j+1}})
$$
and
$$
W^{j+1}=\overline{T_{Z^{j+1\,\prime}}\cap W^j_{\ell_j}}=
\overline{T_{Z^{j+1\,\prime}}\cap V_{k_j}}
$$
is the induced directed structure.

\plainitem{(3.6.5)} for all $j\ge 0$, we have $Z^{j+1}\subsetneq 
Z^j_{\ell_j}$ but $\pi_{k_{j+1},k_{j+1}-1}(Z^{j+1})=Z^j_{\ell_j-1}$.\par}
\smallskip

\noindent
For $j=0$, we have nothing to do by our hypotheses. Assume that $(Z^j,W^j)$
has been constructed. By Proposition 2.4, we get an algebraic subset 
$\smash{\Sigma\subsetneq Z^j_\ell}$
in some stage of the Semple tower $(Z^j_\ell)$ of $Z^j$ such that
every entire curve $f:(\bC,T_\bC)\to(X,V)$ satisfying
$f_{[k_j]}(\bC)\subset Z^j$ also satisfies 
$f_{[k_j+\ell]}(\bC)\subset\Sigma$. By definition, this implies
the first inclusion in the sequence
$$
Z^j=\IEL_{X,V}(Z^j)\subset\pi_{k_j+\ell,k_j}(\IEL_{X,V}(\Sigma))\subset
\pi_{k_j+\ell,k_j}(\Sigma)\subset Z^j
$$
(the other ones being obvious), so we have in fact an equality throughout.
Let $(S_\alpha)$ be the irreducible
components of $\IEL_{X,V}(\Sigma)$. We have $\IEL_{X,V}(S_\alpha)=S_\alpha$ and
one of the components $S_\alpha$ must already satisfy
$\pi_{k_j+\ell,k_j}(S_\alpha)=Z^j=Z^j_0$. We take $\ell_j\in[1,\ell]$ to be 
the smallest order such that $Z^{j+1}:=\pi_{k_j+\ell,k_j+\ell_j}(S_\alpha)
\subsetneq\smash{Z^j_{\ell_j}}$, and set $k_{j+1}=k_j+\ell_j>k_j$. 
By definition of $\ell_j$, we have $\smash{\pi_{k_{j+1},k_{j+1}-1}(Z^{j+1})=
Z^j_{\ell_j-1}}$, otherwise $\ell_j$ would not be minimal.
The fact that $\IEL_{X,V}(S_\alpha)=S_\alpha$ immediately implies
$\IEL_{X,V}(Z^{j+1})=Z^{j+1}$. Also $Z^{j+1}$ cannot
be contained in the vertical divisor $\smash{D_{k_{j+1}}}$. In fact
no irreducible algebraic set $Z$ such that $\IEL_{X,V}(Z)=Z$ can be 
contained in a vertical divisor $D_k$, because $\pi_{k,k-2}(D_k)$ 
corresponds to stationary jets in $X_{k-2}\,$; as every non constant 
curve $f$ has non stationary points, its $k$-jet $f_{[k]}$ cannot
be entirely contained in $D_k$. Finally, the induced directed structure
$(Z^{j+1},W^{j+1})$ must be of general type modulo $X_{k_{j+1}}\to X$,
by the assumption of the theorem and the fact that $k_{j+1}>k_0$.
The inductive procedure is therefore complete.

By Observation 2.2, we have
$$
\rank W^j<\rank W^{j-1}<\ldots<\rank W^1<\rank W^0=\rank V.
$$
After a sufficient number of iterations we reach $\rank W^j=1$. In
this situation the Semple tower of $Z^j$ is trivial,
$\smash{K_{W^j}=W^{j\,*}\otimes \overline{\cJ}_{W^j}}$ is big, and
Proposition~2.4 produces a divisor $\Sigma\subsetneq Z^j_\ell=Z^j$
containing all jets of entire curves with $\smash{f_{[k_j]}(\bC)\subset
Z^j}$.  This
contradicts the fact that $\IEL_{X,V}(Z^j)=Z^j$. We have reached a contradiction,
and Theorem~3.6 is thus proved.\qed

\claim 3.7. Remark|\rm As it proceeds by contradiction, the
proof is unfortunately non constructive -- especially it does not give
any information on the degree of the locus $Y\subsetneq X_{k_0}$ whose
existence is asserted. On the other hand, and this is a bit
surprising, the conclusion is obtained even though the conditions to
be checked do not involve cutting down the dimensions of the base loci
of jet differentials; in fact, the contradiction is obtained even though
the integers $k_j$ may increase and $\dim Z^j$ may become very large.
\endclaim

\noindent
The special case $k_0=0$ of Theorem 3.6 yields the following

\claim 3.8. Partial solution to the generalized GGL conjecture|Let 
$(X,V)$ be a directed pair that is strongly of general type.
Then the Green-Griffiths-Lang conjecture holds true for $(X,V)$, namely
$\ECL(X,V)\subsetneq X$, in other words
there exists a proper algebraic variety $Y\subsetneq X$ such that
every non constant holomorphic curve $f:\bC\to X$ tangent to $V$
satisfies $f(\bC)\subset Y$.
\endclaim

\claim 3.9. Remark|\rm The condition that $(X,V)$ is strongly of 
general type seems to be related to some sort of stability condition.
We are unsure what is the most appropriate definition, but here is one
that makes sense. Fix an ample divisor $A$ on~$X$. For every 
irreducible subvariety $Z\subset X_k$ that projects onto $X_{k-1}$ for
$k\ge 1$, and $Z=X=X_0$ for $k=0$, we define the slope $\mu_A(Z,W)$
of the corresponding directed variety $(Z,W)$ to be
$$
\mu_A(Z,W)=\frac{\inf\lambda}{\rank W},
$$
where $\lambda$ runs over all rational numbers such that there exists
$m\in\bQ_+$ for which 
$$
K_W\otimes\big(\cO_{X_k}(m)\otimes\pi_{k,0}^*\cO(\lambda A)
\big)_{|Z}\quad\hbox{is big on $Z$}
$$
(again, we assume here that $Z\not\subset D_k$ for $k\ge 2$). Notice that
$(X,V)$ is of general type if and only if $\mu_A(X,V)<0$, and that
$\mu_A(Z,W)=-\infty$ if $\cO_{X_k}(1)_{|A}$ is big. Also, the proof
of Lemma~1.11 shows that 
$$
\mu_A(X_k,V_k)\le \mu_A(X_{k-1},V_{k-1})\le\ldots\le
\mu_A(X,V)\quad\hbox{for all~$k$}
$$
(with $\mu_A(X_k,V_k)=-\infty$ for $k\ge k_0\gg 1$ if $(X,V)$ is of
general type). We say that $(X,V)$ is {\it $A$-jet-stable} (resp.\ 
{\it $A$-jet-semi-stable})
if $\mu_A(Z,W)<\mu_A(X,V)$ (resp.\ $\mu_A(Z,W)\le\mu_A(X,V)$) for
all $Z\subsetneq X_k$ as above. It is then clear that if
$(X,V)$ is of general type and $A$-jet-semi-stable, then it is strongly
of general type in the sense of Definition~3.1. It would be useful 
to have a better understanding of this condition of stability 
(or any other one that would have better properties).\qed
\endclaim

\claim 3.10. Example: case of surfaces|\rm Assume that $X$ is a 
minimal complex surface of general type and $V=T_X$ (absolute case). 
Then $K_X$ is nef and big and the Chern classes of $X$ satisfy
$c_1\le 0$ ($-c_1$ is big and nef) and $c_2\ge 0$. The Semple jet-bundles 
$X_k$ form here a tower of $\bP^1$-bundles and $\dim X_k=k+2$.
Since $\det V^*=K_X$ is big, the strong general type assumption of 
3.6 and 3.8 need only be checked for irreducible hypersurfaces
$Z\subset X_k$ distinct from $D_k$ that project onto $X_{k-1}$, of 
relative degree~$m$. The projection $\pi_{k,k-1}:Z\to X_{k-1}$ is a 
ramified $m:1$ cover. Putting $\cO_{X_k}(Z)\simeq\cO_{X_k}({\bf a})\otimes
\pi_{k,0}(B)$, $B\in\Pic(X)$, we can apply (3.4.1) to get an inclusion
$$
K_W\supset
\big(\cO_{X_k}({\bf a}-{\bf1})\otimes\pi_{k,0}^*(B\otimes K_X)\big)_{|Z}
\otimes\cJ_S,\qquad
{\bf a}\in\bZ^k,~~a_k=m.
$$
Let us assume $k=1$ and $S=\emptyset$ to make things even simpler, and let us
perform numerical calculations in the cohomology ring
$$H^\bullet(X_1,\bZ)=H^\bullet(X)[u]/(u^2+c_1u+c_2),\qquad u=c_1(O_{X_1}(1))
$$
(cf.\ [DEG00, Section~2] for similar calculations and more details). We have
$$
Z\equiv mu+b\quad\hbox{where}\quad b=c_1(B)\quad\hbox{and}\quad
K_W\equiv(m-1)u+b-c_1.$$
We are allowed here to add to $K_W$ an arbitrary multiple 
$\cO_{X_1}(p)$, $p\ge 0$, which we rather write $p=mt+1-m$, $t\ge 1-1/m$.
An evaluation of the Euler-Poincar\'e characteristic of $K_W+\cO_{X_1}(p)_{|Z}$ requires computing the intersection number
$$
\eqalign{
\big(K_W+\cO_{X_1}(p)_{|Z}\big)^2\cdot Z&=\big(mt\,u+b-c_1\big)^2(mu+b)\cr
&=m^2t^2\big(m(c_1^2-c_2)-bc_1\big)
+2mt(b-mc_1)(b-c_1)+m(b-c_1)^2,\cr}
$$
taking into account  that $u^3\cdot X_1=c_1^2-c_2$. In case $S\ne\emptyset$,
there is an additional (negative) contribution from the ideal $\cJ_S$
which is $O(t)$ since $S$ is at most a curve.
In any case, for $t\gg 1$, the leading term in the expansion is 
$m^2t^2(m(c_1^2-c_2)-bc_1)$ and the other terms are negligible with respect
to $t^2$, including the one coming from $S$. We know that $T_X$ is 
semistable with 
respect to $c_1(K_X)=-c_1\ge 0$. Multiplication by the section $\sigma$ yields a
morphism $\pi_{1,0}^*\cO_X(-B)\to\cO_{X_1}(m)$, hence by direct image,
a morphism $\cO_X(-B)\to S^mT^*_X$. Evaluating slopes against $K_X$
(a big nef class), the semistability condition implies 
$bc_1\le \frac{m}{2}c_1^2$, and our leading term is bigger that 
$m^3t^2(\frac{1}{2}c_1^2-c_2)$. We get a positive anwer in the well-known 
case where $c_1^2>2c_2$, corresponding to $T_X$ being almost
ample. Analyzing positivity for the full range of values $(k,m,t)$ and
of singular sets $S$ seems an unsurmountable task at this point; in general,
calculations made in [DEG00] and [McQ99] indicate that the Chern class
and semistability conditions become less demanding for higher order
jets (e.g.\ $c_1^2>c_2$ is enough for $Z\subset X_2$, and 
$c_1^2>\frac{9}{13}c_2$ suffices for $Z\subset X_3$). When $\rank V=1$,
major gains come from the use of Ahlfors currents in combination
with McQuillan's tautological inequalities [McQ98]. We therefore hope for a 
substantial strengthening of the above sufficient conditions, and 
a better understanding of the stability issues, possibly in combination
with a use of  Ahlfors currents and tautological inequalities. In the case
of surfaces, an application of Theorem 3.6 for $k_0=1$ and an analysis
of the behaviour of rank $1$ (multi-)foliations on the surface~$X$
(with the crucial use of [McQ98]) was the main argument used in [DEG00]
to prove the hyperbolicity of very general surfaces of degree $d\ge 21$
in $\bP^3$. For these surfaces, one has $c_1^2<c_2$ and $c_1^2/c_2\to 1$ as
$d\to+\infty$. Applying Theorem 3.6 for higher values $k_0\ge 2$ might 
allow to enlarge the range of tractable surfaces, if the behavior of
rank $1$ (multi)-foliations on $X_{k_0-1}$ can be analyzed independently.
\endclaim

\plainsection{4. Algebraic jet-hyperbolicity implies Kobayashi hyperbolicity}

Let $(X,V)$ be a directed variety, where $X$ is an irreducible projective 
variety; the concept still makes sense when $X$ is singular, by embedding
$(X,V)$ in a projective space $(\bP^N,T_{\bP^N})$ and taking the linear space
$V$ to be an irreducible algebraic subset of $T_{\bP^n}$ that is contained 
in $T_X$ at regular points of~$X$.

For any irreducible algebraic subvariety $Z\subset X$, we get as in
section~2 a directed variety~structure $(Z,W)\subset (X,V)$ by taking
$W=\overline{T_{Z'}\cap V}$ on a sufficiently small Zariski open
set $Z'\subset Z_{\rm reg}$ where the intersection has minimal rank. 
Notice that when $W=0$ there cannot exist entire curves $f:(\bC,T_\bC)\to
(Z,W)$, except possibly those which lie in the algebraic set 
$Z\smallsetminus Z'$, hence this case is easy to deal with by induction on
dimension. Otherwise, we can resolve singularities of $Z$ to get a directed 
variety $(\smash{\widehat Z},\smash{\widehat W})$ where $\smash{\widehat Z}$
is nonsingular and $\rank\smash{\widehat W}\ge 1$.

\claim 4.1. Definition|Let $(X,V)$ be a directed variety. We say that 
{\plainitemindent=7mm
\vskip3pt
\plainitem{\rm (a)} $(X,V)$ is algebraically jet-hyperbolic
if for every irreducible algebraic subvariety $Z_0\subset X$, the induced
directed structure $(Z_0,W_0)$ either satisfies $W_0=0$, or has a 
desingularization
$(\smash{\widehat Z_0},\smash{\widehat W_0})$, 
$\rank\smash{\widehat W_0}\ge 1$, that is strongly of general type.
\vskip3pt
\plainitem{\rm (b)} $(X,V)$ is algebraically fully jet-hyperbolic if
for every $k\ge 0$ and every irreducible algebraic subvariety 
$Z\subset X_k$ that is not contained in the union $\Delta_k$ of vertical
divisors, the induced directed structure $(Z,W)$ either satisfies 
\hbox{$W=0$}, or is of general type modulo $X_k\to X$,
i.e.\ has a desingularization $(\smash{\widehat Z},\smash{\widehat W})$, 
$\mu:\smash{\widehat Z}\to Z$, such that some twisted canonical sheaf
$K_{\widehat W}\otimes\mu^*(\cO_{X_k}({\bf a})_{|Z})$, ${\bf a}\in\bN^k$, 
is big. \vskip0pt}
\endclaim

It is clear that hypothesis 4.1~(b) is stronger than 4.1~(a).
In fact, in 4.1~(a), one first takes an induced directed subvariety
$(Z_0,W_0)\subset(X,V)$, its Semple tower $(Z_k,W_k)$, and the question is
to check whether every induced subvariety $(Z,W)\subset(Z_k,W_k)\subset
(X_k,V_k)$ such that $\pi_{k,k-1}(Z)=Z_{k-1}$ is of general type modulo
$Z_k\to Z_0$. On the other hand, for property 4.1~(b), we have to check
right away all induced structures $(Z,W)\subset(X_k,V_k)$, whatever are
their projections $\pi_{k,\ell}(Z)$. It is unclear to us whether the
resulting concepts are really different. Thanks to Theorem 3.8, a 
very easy induction on the dimension of $X$ implies

\claim 4.2. Theorem|Let $(X,V)$ be an irreducible projective directed 
variety that is algebraically jet-hyperbolic in the sense of the above 
definition. Then $(X,V)$ is Brody $($or Kobayashi$\,)$ hyperbolic, i.e.\ 
$\ECL(X,V)=\emptyset$.
\endclaim

\proof By Theorem 3.8, we have $Y:=\ECL(X,V)\subsetneq X$. If $Y\ne\emptyset$, apply induction on dimension to each of the irreducible components
$X'_j$ of $Y$ and to the induced directed structures $(X'_j,V'_j)$ to get
$\ECL(X,V)\subset\bigcup\ECL(X'_j,V'_j)\subsetneq \bigcup X'_j=Y$, a 
contradiction.\qed

\plainsection{5. Proof of the Kobayashi conjecture on generic hyperbolicity}

We start with a general situation, and then restrict ourselves to the
special case of complete intersections in projective space. Consider a
smooth deformation $\pi:\cX\to S$ of complex projective manifolds,
i.e.\ a proper algebraic submersion over a quasi-projective algebraic
manifold $S$ such that the fibers are nonsingular. By a ``very general
fiber'', we mean here a fiber $X_t=\pi^{-1}(t)$ over a point $t$ taken
in the complement $S\smallsetminus\bigcup S_\nu$ of a countable union
of algebraic subsets $S_\nu\subsetneq S$. We are only interested in
the very general fiber and can therefore restrict ourselves to the case
where $S$ is affine after replacing $S$ with a suitable Zariski open
subset~$S^0\subset S$. Ample vector bundles over the total space 
$\cX$ are then the same as vector bundles that are relatively ample 
over~$S$, as one can see immediately by the direct image theorem and
the fact that every locally free sheaf on an affine variety is very ample.

\claim 5.1. Theorem|Let $\pi:\cX\to S$ be a deformation of complex
projective nonsingular varieties $X_t=\pi^{-1}(t)$ over a smooth
quasi-projective irreducible base $S$. Let $n=\dim X_t$ be the
relative dimension and let $N=\dim S$.  Assume that for all
$q=N+1,\ldots,N+n$, the exterior power $\Lambda^qT^*_{\cX}$ is a
relatively ample vector bundle over~$S$.  Then the very general fiber
$X_t$ is algebraically jet-hyperbolic, and thus Kobayashi
hyperbolic.
\endclaim

\proof By taking the relative directed structure $\cV=T_{\cX/S}=
\Ker(d\pi:T_\cX\to\pi^*T_S)$ on~$\cX$,
one constructs a ``relative'' Semple tower $(\cX_k,\cV_k)$ over $\cX$.
It specializes to the absolute Semple tower $X_{t,k}$ of 
$X_t=\pi^{-1}(t)\subset\cX$ when one takes 
$X_{t,k}=\pi_{k,0}^{-1}(X_t)\subset \cX_k$ via the natural projection
$\pi_{k,0}:\cX_k\to\cX_0=\cX$. By construction $\cV_0=\cV=T_{\cX/S}$ and all
$\cV_k$ have rank $n$. Let $(\cX^a_k,\cV^a_k)$ be the absolute Semple 
tower of $\cX$, so that $\cX^a_0=\cX$ and $\cV^a_0=T_{\cX}$,
and let $\smash{\widetilde\cV}_k$ be the
restriction of the vector bundle $\cV^a_k$ to $\cX_k\subset\cX^a_k$, so that
$\rank\smash{\widetilde\cV}_k=\rank\cV^a_k=N+n$.
For~every $k\ge 0$, we claim that there is an exact sequence of 
vector bundles 
$$0\to \cV_k\to \widetilde\cV_k\to \cS_k\to 0,\qquad
\cS_k\simeq(\pi\circ\pi_{k,0})^*T_S\otimes\cO_{\cX_k}({\bf 1})\quad
\hbox{over $\cX_k$},
\leqno(5.1.1)$$
where ${\bf 1}=(1,\ldots,1)\in\bN^k$, $\rank\cV_k=n$, and
$\rank\widetilde\cV_k=\rank\cV^a_k=N+n=\dim\cX$. Since 
$\smash{\widetilde\cV_0}=\cV^a_0=T_\cX$ and $\cV_0=T_{\cX/S}$, this is true 
by definition for $k=0$, with $\cS_0=\pi^*T_S$ and 
$\cO_{\cX_0}({\bf 1})=\cO_\cX$.
In general, there is a well defined injection of bundles
$\cV_k\to \smash{\widetilde\cV_k}$, the quotient is of rank~$N$, and 
we simply put $\cS_k=\smash{\widetilde\cV_k}/\cV_k$ by definition.
The relative (resp.\ absolute) Semple tower of $\pi_{k,\ell}:\cX_k\to\cX_\ell$ 
(resp.\ $\pi^a_{k,\ell}:\cX^a_k\to\cX^a_\ell$) yields exact sequences
$$\leqalignno{
&0\longrightarrow\cG_k\longrightarrow \cV_k
\mathop{\longrightarrow}\limits^{(\pi_{k,k-1})_*}\cO_{\cX_k}(-1)
\longrightarrow 0,\qquad\cG_k:=T_{\cX_k/\cX_{k-1}},
\raise-4pt\hbox{\strut}&(5.1.2)\cr
&0\longrightarrow \cO_{\cX_k}\longrightarrow(\pi_{k,k-1})^*\cV_{k-1}\otimes
\cO_{\cX_k}(1)\longrightarrow\cG_k\longrightarrow 0,&(5.1.3)\cr
&0\longrightarrow\cG^a_k\longrightarrow \cV^a_k
\mathop{\longrightarrow}\limits^{(\pi^a_{k,k-1})_*}\cO_{\cX^a_k}(-1)
\longrightarrow 0,\qquad\cG^a_k:=T_{\cX^a_k/\cX^a_{k-1}},
\raise-4pt\hbox{\strut}&(5.1.2^a)\cr
&0\longrightarrow \cO_{\cX^a_k}\longrightarrow(\pi^a_{k,k-1})^*\cV^a_{k-1}\otimes
\cO_{\cX^a_k}(1)\longrightarrow\cG^a_k\longrightarrow 0.&(5.1.3^a)\cr}$$
By restricting the absolute ones to $\cX_k\subset\cX^a_k$ and denoting
$\smash{\widetilde\cG}_k:=\cG^a_{k|\cX_k}$, we get exact sequences
$$\leqalignno{
&0\longrightarrow\widetilde\cG_k\longrightarrow \widetilde\cV_k
\mathop{\longrightarrow}\limits^{(\pi_{k,k-1})_*}\cO_{\cX_k}(-1)
\longrightarrow 0,\raise-4pt\hbox{\strut}&(5.1.2^\sim)\cr
&0\longrightarrow \cO_{\cX_k}\longrightarrow(\pi_{k,k-1})^*\widetilde
\cV_{k-1}\otimes\cO_{\cX_k}(1)
\longrightarrow\widetilde\cG_k\longrightarrow 0&(5.1.3^\sim)\cr}$$
There is an inclusion morphism of $(5.1.i)$ into $(5.1.i^\sim)$, $i=2,3$,
and by taking cokernels, we see that 
$$
\cS_k:=\widetilde\cV_k/\cV_k\mathop{=\,}\limits_{(*)}\widetilde\cG_k/\cG_k
\mathop{=\,}\limits_{(**)}(\pi_{k,k-1})^*\cS_{k-1}\otimes \cO_{\cX_k}(1)
$$
where $(*)$ comes from $(5.1.2^\sim)$ and $(**)$ from $(5.1.3^\sim)$.
This induction formula for $\cS_k$ completes the proof of (5.1.1).
If we take the dual exact sequences, we get
$$
\leqalignno{
&0\longrightarrow \cO_{\cX_k}(1)\longrightarrow \widetilde\cV_k^*
\longrightarrow\widetilde\cG^*_k\longrightarrow 0,\raise-6pt\hbox{\strut}&(5.1.2^*)\cr
&0\longrightarrow \widetilde\cG^*_k\longrightarrow(\pi_{k,k-1})^*\widetilde\cV^*_{k-1}\otimes
\cO_{\cX_k}(-1)\longrightarrow \cO_{ \cX_k}\longrightarrow 0,&(5.1.3^*)\cr}
$$
and the $q$-th (resp.\ $q'$-th) exterior power of these yield
$$
\leqalignno{
&0\longrightarrow\Lambda^{q-1}\widetilde\cG^*_k\otimes\cO_{\cX_k}(1)\longrightarrow 
\Lambda^q\widetilde\cV_k^*
\longrightarrow\Lambda^q\widetilde\cG^*_k\longrightarrow 0,\raise-6pt\hbox{\strut}
&(5.1.4)\cr
&0\longrightarrow \Lambda^{q'}\widetilde\cG^*_k\longrightarrow(\pi_{k,k-1})^*\Lambda^{q'}
\widetilde\cV^*_{k-1}
\otimes\cO_{\cX_k}(-q')\longrightarrow \Lambda^{q'-1}\widetilde\cG^*_k\longrightarrow 0.
&(5.1.5)\cr}
$$
In a next step, we will need local vector fields and their liftings to the 
absolute and relative Semple towers to justify certain delicate arguments
about multiplier ideals.

\claim 5.1.6. Lemma|
\vskip3pt
{\plainitemindent=7mm
\plainitem{\rm (a)} Every local holomorphic vector field $\zeta$ on an open
set $U\subset \cX$ has a natural lifting $\zeta^{(k)}$ to the open set
$(\pi^a_{k,0})^{-1}(U)$ in the total space of the absolute Semple 
bundle $\cX_k^a$. 
\vskip3pt
\plainitem{\rm (b)} 
In particular, if one assumes that $\zeta$ is 
in~$H^0(U,T_{\cX/S})\subset H^0(U,T_{\cX})$, the lifted flow leaves 
invariant the relative Semple tower
$\cX_k\subset\cX_k^a$ in $\pi_{k,0}^{-1}(U)$, thus 
$\smash{\zeta^{(k)}_{|\cX_k}}$ is tangent to~$\cV_k$. 
\vskip3pt
\plainitem{\rm (c)} 
Every local holomorphic vector field $\tau\in H^0(\Omega,T_S)$ on $S$
has a $($nonunique$\,)$ lifting $\widetilde\tau\in H^0(U,T_\cX)$ on a neighborhood
$U$ of every point $x\in\pi^{-1}(\Omega)$. Once $\widetilde\tau$ is chosen, there
is a $($unique$\,)$ lifting $\widetilde\tau^{(k)}$ of $\widetilde\tau$ to the open set
$(\pi^a_{k,0})^{-1}(U)$ in the total space of the absolute
Semple bundle $\cX^a_k$, and the flow of $\smash{\widetilde\tau^{(k)}}$ induces
a local biholomorphism $(\cX_k,\cV_k)\to(\cX_k,\cV_k)$ of the relative Semple
tower. These local flows preserve the exact sequences $(5.1.1-5.1.5)\,;$
moreover $\smash{\widetilde\tau^{(k)}_{|\cX_k}}$ is tangent to~$\cV^a_k$ $($but is
not tangent to $\cV_k\,)$.
\vskip0pt}
\endclaim

\proof (a) Every local holomorphic vector field $\zeta$ on $\cX$ generates a
flow of local biholomorphisms on open subsets of $\cX$, and 
we can apply the fonctoriality property~1.10 to lift it to a flow 
on $(\cX_k^a,\cV_k^a)$. The differentiation of the lifted flow gives 
back what we define to be the lifted vector field $\zeta^{(k)}$ 
on~$\cX_k^a$.\vskip3pt

\noindent(b) If $\zeta$ lies in $\cV_0=T_{\cX/S}$, the 
lifted flow acts similarly on the relative tower $(\cX_k,\cV_k)$ since
the fibers $X_{t,k}$ over $X_t=\pi^{-1}(t)$ are preserved. Therefore
$\smash{\zeta^{(k)}_{|\cX_k}}$ is tangent to $\cV_k$, since the relative
Semple tower is the absolute Semple tower of the fibers~$X_t$.
\vskip3pt

\noindent(c) Since $\cX\to S$ 
is a holomorphic submersion, $\cX$~is locally holomorphically
trivialized as a product $S\times X_t\,$; if $\tau$ is a local vector 
field on~$S$, it can thus be lifted (possibly not uniquely) as a 
local vector field $\widetilde\tau$ on $\cX$. The resulting flow of 
$\widetilde\tau$ on $\cX$ commutes with the flow of $\sigma$ on $S$, 
i.e.\ acts by (local) biholomorphisms $(\cX,T_{\cX/S})\to(\cX,T_{\cX/S})$ of
directed varieties. From there we conclude, again by
fonctoriality, that the lifting of the flow to $\cX_k^a$
preserves $(\cX_k,\cV_k)\subset (\cX^a_k,\cV^a_k)$, in particular 
the exact sequences are preserved. Moreover the differential
$\smash{\widetilde\tau^{(k)}}$ is tangent to~$\cV^a_k$ (but, already for $k=0$,
it is not tangent to~$\cV_k$).\qed\vskip5pt

Now, let $\cZ\subset\cX_k$ be an irreducible algebraic subvariety
of $\cX_k$ that is not contained in the union $\Delta_k$ of vertical
divisors and projects surjectively onto~$S$, and let 
$(\cZ,\cW)\subset(\cX_k,\cV_k)$ be the induced directed structure,
Our ultimate goal is to show that the generic fibers $(Z_t,W_t)$ are
of general type modulo~$X_{t,k}\to X_t$. We assume $r:=\rank\cW\ge 1$,
otherwise there is nothing to do. We first need to extend $\cZ$ 
within the absolute Semple tower~$\cX^a_k$, ``horizontally'' with
respect to the projection $\cX^a_k\to S$. Later, we will have to care 
that the extension $\cZ^a$ is made in an ``equisingular way''. This
procedure is what replaces here the use of oblique vector fields
(introduced by Siu [Siu02, 04], and employed later by [Pau08] and
[DMR10]). The main advantage is that ``equisingular horizontal extensions''
do not introduce additional poles.

\claim 5.1.7. Lemma|There exists an irreducible algebraic variety
$\cZ^a\subset\cX^a_k$ satisfying the following properties with respect
to the projections $\cZ_\ell=\pi_{k,\ell}(\cZ)$ and 
$\cZ^a_\ell=\pi_{k,\ell}(\cZ^a)$, $\ell=0,1,\ldots,k\,:$
{\plainitemindent=7mm
\vskip3pt
\plainitem{\rm(a)} $\cZ$ is one of the irreducible components of
$\cZ^a\cap\cX_k$ and likewise, $\cZ_\ell$ is one of the irreducible 
components of $\cZ^a_\ell\cap\cX_\ell$ for all $\ell\,;$
\vskip3pt
\plainitem{\rm(b)} the intersection $\cZ^a_\ell\cap\cX_\ell$ is smooth and
transverse at the generic point of $\cZ_\ell\,;$
\vskip3pt
\plainitem{\rm(c)} the codimension of $\cZ^a_\ell$ in $\cX^a_\ell$ is
equal to the codimension $p_\ell$ of $\cZ_\ell$ in $\cX_\ell$ for all~$\ell\,;$
\vskip3pt
\plainitem{\rm(d)} if $r_\ell$ is the $($generic$\,)$ rank 
of $\cW_\ell=T_{\cZ_\ell}\cap\cV_\ell$, then the rank of
$\cW^a_\ell=T_{\cZ^a_\ell}\cap\cV^a_\ell$ is $N+r_\ell\,;$
\vskip3pt
\plainitem{\rm(e)} we have $r_\ell=n-p_\ell$, $p_0\le p_1\le\ldots\le p_k$ and
$r_0\ge r_1\ge\ldots\ge r_k=r\ge 1$.
\vskip0pt}
\endclaim

\proof For $\ell=0$, no extension is needed as $\cX^a_0=\cX_0=\cX$, we simply 
take $\cZ^a_0=\cZ_0=\pi_{k,0}(\cZ)$ (and thus we are done if $k=0$). 
For $k\ge 1$, we construct $\cZ^a_\ell$ inductively, assuming that
$\cZ^a_{\ell-1}$ has already been constructed, $\ell\ge 1$. Since
$\cZ_{\ell-1}=\pi_{k,k-1}(\cZ_\ell)$ and $\pi_{\ell,\ell-1}:\cX_\ell\to\cX_{\ell-1}$ is
a fibration, it is clear that $p_{\ell-1}\le p_\ell$. Also, $\cZ_\ell$ is
contained in $P(\cW_{\ell-1})\subset P(\cV_{\ell-1})=\cX_\ell$. The codimension
of $P(\cW_{\ell-1})$ in $\cX_\ell$ is $p_{\ell-1}=n-r_{\ell-1}$, hence the
codimension of $\cZ_\ell$ in $P(\cW_{\ell-1})$ is $p_\ell-p_{\ell-1}$.
Now, $\cZ_\ell$ is the zero locus of a family of sections $s_j$ of some
very ample line bundle $\cO_{\cX_\ell}({\bf m}_\ell)$ on $\cX_\ell$, 
${\bf m}_\ell\in\bN^\ell$, such that the differentials $ds_j$ are 
independent at all nonsingular points of~$\cZ_\ell$. If we take
${\bf m}_\ell$ large enough, those sections $s_j$ extend as sections $s^a_j$
of $\cO_{\cX^a_\ell}({\bf m}_\ell)$, and we can pick $p_\ell-p_{\ell-1}$ linear
combinations $\sigma_j$ of them so that $P(\cW_{\ell-1})\cap\{\sigma_j=0\}$ and
$P(\cW^a_{\ell-1})\cap\{\sigma^a_j=0\}$ are generically transverse intersections
of pure dimension, the dimension of $P(\cW_{\ell-1})\cap\{\sigma_j=0\}$ being
equal to that of~$\cZ_\ell$. We take $\cZ^a_\ell$ to be the 
irreducible component of $P(\cW^a_{\ell-1})\cap\{s^a_j=0\}$ that 
contains~$\cZ_\ell$ (in general, the intersection will not be irreducible,
as its degree may be very large). Properties (a), (b), (c) are then satisfied
by construction, and $\cW_\ell$ (resp.\ $\cW^a_\ell$) is obtained from
the lifting of $\cW_{\ell-1}$ to $P(\cW_{\ell-1})$ (resp.\ of
$\cW^a_{\ell-1}$ to $P(\cW^a_{\ell-1})$) by cutting the corresponding lifted
directed structure by the generically independent linear equations 
$d\sigma_j=0$ (resp.\ $d\sigma^a_j=0$). Therefore we see by induction
that the rank of $\cW_\ell$ (resp.~$\cW^a_\ell$) is 
$$r_\ell=r_{\ell-1}-(p_\ell-p_{\ell-1})=n-p_\ell,\qquad\hbox{resp.}\quad
N+r_{\ell-1}-(p_\ell-p_{\ell-1})=N+r_\ell.$$
Properties (d), (e) follow, and Lemma~5.1.7 is proved.\qed
\vskip4pt

To simplify notation, we let $(\cZ^a,\cW^a)=(\cZ^a_k,\cW^a_k)$ be
the induced directed structure at the top level~$k$ (the lower levels will
no longer be needed in the sequel), i.e.\ $\cW^a=T_{\cZ^a}\cap\cV^a_k$ at 
the generic point. We take $\smash{\widetilde\cW}$ to be the irreducible 
component of $\smash{\cW^a_{|\cZ}}$ that contains $\cW$ (recall that 
$\cZ\subset\cZ^a\cap\cX_k$ have the same dimension, but the intersection
need not be irreducible). We claim that there is an exact sequence
$$
0\longrightarrow \cW\longrightarrow\smash{\widetilde\cW}
\longrightarrow(\pi\circ\pi_{k,0})^*T_S\otimes\cO_{\cX_k}({\bf 1})
\longrightarrow 0\leqno(5.1.8)
$$
at a generic point of~$\cZ$. In fact, at a generic point, the proof of
Lemma~5.1.7 shows that
$$\cW=\cV_k\cap\{d\sigma_j=0\},\quad
\cW^a=\cV^a_k\cap\{d\sigma^a_j=0\},\quad
\widetilde\cW=\widetilde\cV_k\cap\{d\sigma^a_j=0\}$$
have the same codimension  in $\cV_k$, $\cV^a_k$,
$\smash{\widetilde\cV}_k$ respectively, hence we have
$\smash{\widetilde\cW}/\cW\simeq\smash{\widetilde\cV}_k/\cV_k\simeq\cS_k$,
and the conclusion follows by restricting the exact sequence~(5.1.1).
We get by definition 
a non trivial morphism over~$\cZ$ induced by the natural inclusion
$\smash{\widetilde\cW\subset\widetilde\cV_k}$
$$\Lambda^q\widetilde\cV^*_{k\,|\cZ}\longrightarrow K_{\widetilde\cW},\qquad
q=N+r=\rank\widetilde\cW.$$
One should notice that $\widetilde\cV_k$ and $\widetilde\cG_k$ are genuine vector bundles
without singularities, hence the above morphism actually has its image 
{\it contained in} $K_{\widetilde\cW}=K_{\cW^a|\cZ}$ even when one takes into 
account the relevant multiplier ideal sheaf that 
defines~$K_{\cW^a}$ (``monotonicity principle''); however, there
will be more delicate singularity issues later on. We conclude by 
$(5.1.4)$ that either we have a non trivial morphism
$$\Lambda^{q-1}\widetilde\cG^*_k\otimes\cO_{\cX_k}(1)_{|\cZ}\longrightarrow K_{\widetilde\cW}$$
or (if the above vanishes) a non trivial morphism
$$\Lambda^q\widetilde\cG^*_{k\,|\cZ}\longrightarrow K_{\widetilde\cW}.$$
By $(5.1.5)$ with $q'=q$ or $q'=q+1$, we infer that we have a non trivial
morphism
$$
(\pi_{k,k-1})^*\Lambda^q\widetilde\cV^*_{k-1}\otimes\cO_{\cX_k}(-q+1)_{|\cZ}
\longrightarrow K_{\widetilde\cW}
$$
or a non trivial morphism
$$
(\pi_{k,k-1})^*\Lambda^{q+1}\widetilde\cV^*_{k-1}\otimes\cO_{\cX_k}(-q-1)_{|\cZ}
\longrightarrow K_{\widetilde\cW}.
$$
Proceeding inductively with the lower stages and getting down to
$\widetilde\cV_0=T_{\cX}$, we conclude that there exists an integer
$q'\ge q=\rank\smash{\widetilde\cW}$, a weight 
${\bf a}=(a_1,\ldots,a_k)\in\bN^k$, $a_j\ge q-1$, and a non trivial 
morphism
$$
(\pi_{k,0})^*\Lambda^{q'}T^*_{\cX\,|\cZ}\to K_{\widetilde\cW}
\otimes\cO_{\cX_k}({\bf a})_{|\cZ}.\leqno(5.1.9)
$$
Assume that $q\ge N+1$ (i.e.\ $\rank\cW\ge 1$). By our assumption 
(assuming $S$ affine here), $\Lambda^{q'}T^*_{\cX}$ is ample 
over~$\cX$, thus, by twisting with a certain relatively ample line 
bundle $\cO_{\cX_k}(\varepsilon{\bf c})$ with respect to $\pi_{k,0}$, we see that
$(\pi_{k,0})^*\Lambda^{q'}T^*_{\cX}\otimes\cO_{\cX_k}(\varepsilon{\bf c})$ is ample
over $\cX_k$ for $0<\varepsilon\ll 1$. From this, we infer that there exists 
a weight ${\bf b}\in\bQ_+^k$, $b_j>q-1$, such that 
$K_{\widetilde\cW}\otimes\cO_{\cX_k}({\bf b})_{|\cZ}$
is big over $\cZ$. In fact, these arguments could have been given instead
for $(\cZ^a,\cW^a)$, and we would have obtained in the same manner the 
existence of a non trivial morphism
$$
(\pi^a_{k,0})^*\Lambda^{q'}T^*_{\cX\,|\cZ^a}\to K_{\cW^a}
\otimes\cO_{\cX^a_k}({\bf a})_{|\cZ},\leqno(5.1.9^a)
$$
in other words $(\cZ^a,\cW^a)$ is of general type modulo 
$\cX^a_k\to\cX$. The reasoning is of course much easier for $k=0$,
in that case we simply take $\cZ^a=\cZ$, $\cW^a=T_{\cZ}$, $q'=q$,
and $\Lambda^qT^*_{\cX}$ is ample on~$\cX$, hence on $\cZ$.
\vskip4pt

From this, we are going to conclude by a Hilbert scheme argument
that $(X_t,T_{X_t})$ is algebraically jet-hyperbolic for very general 
$t\in S$. Otherwise, consider the collection of non vertical
irreducible varieties
$Z_t\times\{t\}\subset\cX_k$ such that the induced directed structure $(Z_t,W_t)$ is not of general type modulo $X_{t,k}\to X_t$, with $\rank W_t\ge 1$ and
$t$ running over~$S$. If we fix $k$, the degree $\delta$ of $Z_t$ with
respect to some polarization and the weight ${\bf b}\in\bQ_+^k$ such that
$K_{W_t}\otimes\cO_{\cX_k}({\bf b})_{|Z_t}$ is not big, we get a Zariski closed
set $\cH_{k,\delta,{\bf b}}$ in the Hilbert scheme of $\cX_k$, and so is
$\cH_{k,\delta}=\bigcap_{\bf b}\cH_{k,\delta,{\bf b}}$. We have a natural projection
$p_{k,\delta}:\cH_{k,\delta}\to S$. If $p_{k,\delta}$ were dominant, it would
be possible to find a Zariski open set $S^0\subset S$, a finite unramified
cover $\smash{\widehat S}^0$ of $S^0$ and a branched section 
$\smash{\widehat S}^0\to \cH_{k,\delta}$ of $p_{k,\delta}$. This would give
an algebraic family 
$Z_t\subset X_{t,k}$ for $t\in\smash{\widehat{S}^0}$, 
such that the induced directed structure $(Z_t,W_t)$ is not of general type 
modulo $X_{t,k}\to X_t$, with $\rank W_t\ge 1$. In order to avoid finite
covers of the base, we apply a base change 
$\smash{\widehat S}^0\to S$ and consider the resulting deformation
$\smash{\widehat\cX\to\widehat S^0}$, which we still denote
$\cX\to S$ to simplify notation (so that we just have 
$\smash{\widehat{S}^0}=S$ in the new setting). In this way,
we obtain a directed subvariety $(\cZ,\cW)$ of $(\cX_k,\cV_k)$, and we
extend it horizontally as a subvariety $(\cZ^a,\cW^a)\subset (\cX^a_k,\cV^a_k)$
satisfying the following properties:
$$\leqalignno{
\kern30pt&\cZ:=\overline{\bigcup_{t\in S}Z_t}\subset \cX_k\subset\cX^a_k,\quad
\cW:= \overline{T_{\cZ_{\rm reg}}\cap \cV_k}
\subset\cV_{k\,|\cZ}\kern9pt
\hbox{restricts to $W_t$ on $Z_t$,~~$t\in S$},&(5.1.10)\cr
&\hbox{$\cZ$ is one of the irreducible components of $\cZ^a\cap\cX_k$,}
&(5.1.11)\cr
\noalign{\vskip6pt}
&\hbox{if one takes~~$\widetilde\cW=\cW^a_{|\cZ}$
where $\cW^a=\overline{T_{\cZ^a_{\rm reg}}\cap\cV^a_k}$, there is an exact 
sequence}
&(5.1.12)\cr
&\,0\to\cW\to\widetilde\cW\to(\pi\circ\pi_{k,0})^*T_S
\otimes\cO_{\cX_k}({\bf 1})_{|\cZ}\to0\quad\hbox{generically on $\cZ$}
&\cr}
$$
(maybe after shrinking again $S$ to a smaller Zariski open set). The next
idea, which refines the technique used by [Ein88, 91] and [Voi96], is to
consider an embedded desingularization of our spaces to take
care of the singularities and their associated multiplier ideal sheaves.

\claim 5.1.13. Lemma|Let $\mu:\widehat\cZ\to\cZ$ be a desingularization
of $\cZ$, and let $\smash{\widehat Z}_t$ be the fiber over $t\in S$ of 
the projection $\pi\circ\pi_{k,0}\circ\mu:\smash{\widehat\cZ}\to S$. 
In fact, we take $\mu:\smash{\widehat\cX}^a_k\to\cX^a_k$ 
to be a simultaneous embedded resolution of singularities of $\cZ$
and $\cZ^a$ inside the nonsingular space~$\cX^a_k$. Also, by composing
if necessary with further blow-ups, we ask $\mu$ to resolve the 
indeterminacies of the meromorphic maps to Grassmannian bundles
$$\varphi:\cZ\merto\Gr(\cV_k,r),\qquad
\varphi^a:\cZ^a\merto\Gr(\cV^a_k,N+r)$$
associated with the linear spaces $\cW\subset\cV_k\subset T_{\cX_k}$, 
$\cW^a\subset\cV^a_k\subset T_{\cX^a_k}$ of respective ranks
$r$ and~$N+r\,;$ in~other words, we want $\mu$ to provide holomorphic maps
$$\varphi\circ\mu:\widehat\cZ\longrightarrow
\mu^*(\Gr(\cV_k,r)_{|\cZ}),\qquad
\varphi^a\circ\mu:\widehat\cZ^a\longrightarrow
\mu^*(\Gr(\cV^a_k,N+r)_{|\cZ^a}).$$
Under the above assumptions for~$\mu$, let
$(\smash{\widehat\cZ},\smash{\widehat\cW})$,
$(\smash{\widehat\cZ}^a,\smash{\widehat\cW}^a)$ be the pull-backs of 
$(\cZ,\cW)$, $(\cZ^a,\cW^a)$ by~$\mu$, and let
$(\smash{\widehat Z_t},\smash{\widehat W_t})$ be the restriction
of $(\smash{\widehat\cZ}^a,\smash{\widehat\cW}^a)$ to~$\smash{\widehat Z_t}$.
Then, provided that the extension $\cZ\rightsquigarrow\cZ^a$ has been made
in an ``equisingular way'', there is a well defined nontrivial morphism
$$
(\mu^*K_{\widetilde\cW})_{|\widehat Z_t}\longrightarrow K_{\widehat W_t}\otimes
\mu^*\big((\pi\circ\pi_{k,0})^*K_S\otimes\cO_{X_k}(-N\,{\bf 1})\big)_{|\widehat Z_t}
$$
on the generic smooth fiber $\smash{\widehat Z_t}$, $t\in S$, taking into 
account the respective multiplier ideal sheaves of $K_{\widetilde\cW}$
and~$K_{\widehat W_t}$.
\endclaim

\proof We first consider the much easier case~$k=0$ (which does not require
to take an extension $\cZ\rightsquigarrow\cZ^a$). Then 
$\mu:\widehat \cZ\to\cZ\subset\cX\mathop{\to}\limits^\pi S$ is a
fibration over~$S$, we have $\smash{\widetilde W}=T_{\cZ}\subset T_{\cX}$,
and the pull-back of (5.1.12) by $\mu$ reduces to an exact sequence of sheaves
$$
0\to T_{\widehat\cZ/S}\to T_{\widehat\cZ}\to (\pi\circ\mu)^*T_S\to 0.
$$
It restricts to an exact sequence of vector bundles on a neighborhood of
a generic (smooth) fiber $\smash{\widehat Z}_t$, and
$T_{\widehat\cZ/S\,|\widehat Z_t}=T_{\widehat Z_t}=\smash{\widehat W}_t$ by definition.
We then get a composition of morphisms
$$
(\mu^*K_{\widetilde\cW})_{|\widehat Z_t}\longrightarrow K_{\widehat\cZ\,|\widehat Z_t}
\mathop{\longrightarrow}\limits^{\simeq}
K_{\widehat Z_t}\otimes(\pi\circ\mu)^*K_S,
\leqno(5.1.14_0)
$$
which is just the morphism whose existence is asserted in 
the~Lemma. The first arrow is well defined everywhere since by definition
(incorporating multiplier ideals) the canonical sheaf $K_{\widetilde\cW }$
is obtained by restricting smooth sections
of the appropriate exterior power $\Lambda^qT^*_\cX$ to $\cZ$, an operation
followed by pulling-back via a morphism $i_{\cZ}\circ\mu:
\smash{\widehat\cZ}\to\cZ\subset\cX$ between nonsingular varieties.
This completes the case $k=0$.
\vskip4pt
For $k\ge 1$, the situation is more involved: in particular the induced
linear structure $\smash{\widehat\cW}$ on $\smash{\widehat\cZ}$ will
probably remain singular, and $\smash{\widehat W_t}$ may be singular as
well on $\smash{\widehat Z_t}$ even when $\smash{\widehat Z_t}$ is nonsingular.
In any case, (5.1.12) gives an isomorphism
$$
(K_{\widetilde\cW})_{|Z_t}\longrightarrow K_{\cW}\otimes 
(\pi\circ\pi_{k,0})^*K_S\otimes\cO_{X_k}(-N\,{\bf 1})_{\,|Z_t}
$$
at the generic point of $Z_t$ ($t\in S$ being itself generic). By pulling-back
via~$\mu$, we get an isomorphism
$$
(\mu^*K_{\widetilde\cW})_{|\widehat Z_t}\longrightarrow K_{\widehat W_t}\otimes 
\mu^*\big((\pi\circ\pi_{k,0})^*K_S\otimes\cO_{X_k}(-N\,{\bf 1})\big)_{|\widehat Z_t}
\leqno(5.1.14_k)
$$
at the generic point of $\widehat Z_t$. It is obtained from a fibration
between non singular varieties
$$\pi\circ\pi^a_{k,0}\circ\mu:\widehat\cZ\subset\widehat\cX^a_k\to S$$
and we only consider what happens when $\widehat Z_t$ is a nonsingular
fiber. We still have to show that $(5.1.14_k)$ is everywhere defined and 
factorizes through the corresponding multiplier ideals. The hypothesis
that $\mu$ resolves the indeterminacies of the Grassmannian structure maps
of $\cW$ and $\smash{\widetilde\cW}$ implies that we get an exact sequence
of {\it vector bundles}
$$
0\to\mu^*\cW\to\mu^*\widetilde\cW\to\cF\to 0,\quad
\cF\subset \mu^*\cS_{k\,|\widehat\cZ},\quad
\cS_k=(\pi\circ\pi_{k,0})^*T_S\otimes\cO_{\cX_k}({\bf 1})\big)_{|\cZ}
\leqno(5.1.15)
$$
over $\widehat\cZ$ [here $\mu^*\cW$, $\mu^*\widetilde\cW$ mean pull-backs of
{\it linear spaces}, not pull-backs of sheaves, namely, one takes the 
closure of fibers over the open set of regular points]. In fact, we have a 
map $\psi=(\varphi,\varphi^a{}_{|\cZ})$
into the flag bundle $p:F_{r,N+r}(\cV_k,\smash{\widetilde\cV_k})\to\cX_k$ of 
subspaces $L\subset\cV_k$, $L'\subset\smash{\widetilde\cV_k}$ of 
respective dimensions $(r,N+r)$ with $L\subset L'$. 
Let $\cL$ and $\cL'$ be the tautological
vector bundles of rank $r$ and $N+r$ on this flag bundle. By construction
$\psi\circ\mu$ is holomorphic, while $\cW=\psi^*\cL$, 
$\widetilde\cW=\psi^*\cL'$ and $\smash{\widetilde\cV_k}/\cV_k\simeq\cS_k$.
Therefore $\mu^*\cW=(\psi\circ\mu)^*\cL$ and
$\mu^*\smash{\widetilde\cW}=(\psi\circ\mu)^*\cL'$ are genuine vector bundles,
$\mu^*\cW$ is a subbundle of $\mu^*\smash{\widetilde\cW}$ and
$$
\mu^*\smash{\widetilde\cW}/\mu^*\cW=(\psi\circ\mu)^*(\cL'/\cL)
$$
is of rank $N$. On the flag bundle, there is also a natural morphism
$$
\cL'/\cL\to p^*\widetilde\cV_k/p^*\cV_k=p^*\cS_k
$$
between rank $N$ bundles, whose determinant vanishes along a certain 
intrinsically defined divisor~$\Delta$, so that
$\det(\cL'/\cL)=p^*\det\cS_k\otimes\cO(-\Delta)$. In the end, we get
$$\det\cF=\mu^*\Lambda^N\cS_k\otimes(\psi\circ\mu)^*\cO(-\Delta).$$
Now, let $u$ be a local section of $K_{\widetilde\cW}$. By definition 
$u$ is the restriction of a
local holomorphic section $u'$ of $\Lambda^{N+r}T^*_{\cX^a}$ near a
point $x\in\cZ$. In order to construct its image by $(5.1.14_k)$, 
we pick a local generator $\tau_1\wedge\ldots\wedge\tau_N$ of 
$\Lambda^NT_S$ near $s=\pi\circ\pi_{k,0}(x)$. By Lemma~5.1.6~(c),
the local vector fields $\tau_j$ can be lifted as vector fields
$\widetilde\tau_j$ on $\cX$ near $\pi_{k,0}(x)$, and then as vector fields
$\smash{\widetilde\tau_j^{(k)}}$ tangent to $\smash{\cV^a_k\subset T_{\cX^a_k}}$ near~$x$. 
The image of the wedge product 
$\smash{\widetilde\tau_1^{(k)}\wedge\ldots\wedge\widetilde\tau_N^{(k)}{}_{|\cX_k}}$ in
$\Lambda^N\cS_k$ via (5.1.1) is a local generator $\theta$ of 
$$\Lambda^N\cS_k=(\pi\circ\pi_{k,0})^*K_S^{-1}\otimes\cO_{\cX_k}(N\,{\bf 1}),$$
and we want to emphasize here that the construction of the $\tau_j$'s, 
$\widetilde\tau_j$'s, and $\theta$ is made entirely on nonsingular spaces.
Let $\delta$ be a local section of $(\psi\circ\mu)^*\cO(-\Delta)$ on
$\smash{\widehat\cZ}$ whose zero divisor is $(\psi\circ\mu)^*\Delta$. By 
what we have done, $\delta\,\mu^*\theta$ is a local generator
of $\det\cF$, and thanks to (5.1.15), $\delta\,\mu^*\theta$ has a 
local holomorphic lifting $\xi$ with values in 
$\smash{\cO(\mu^*\Lambda^N\widetilde\cW)\subset
\cO(\mu^*\Lambda^N T_{\cX^a_k})}$. We contract 
$\mu^*u'$ with $\xi$ to get a local section $\mu^*u'\cdot\xi$ 
of~$\smash{\Lambda^{r}T^*_{\cX^a_k}}$, and define in this way 
a morphism
$$
\mu^*u\longmapsto (\mu^*u'\cdot\xi)\otimes
(\delta\,\mu^*\theta)^{-1}{}_{|\widehat Z_t}
\leqno(5.1.16)
$$
with values in
$$
K_{\widehat W_t}\otimes 
\mu^*\big((\pi\circ\pi_{k,0})^*K_S\otimes\cO_{X_k}
(-N\,{\bf 1})\big)\otimes\cO((\psi\circ\mu)^*\Delta)_{|\widehat Z_t}.
\leqno(5.1.17)
$$
Observe that $\smash{(\mu^*u'\cdot\theta')_{|\widehat Z_t}}$
is the restriction of a section of the ambient $r$-form bundle 
$$\mu^*\Lambda^{r}T^*_{\cX^a_k}\subset\Lambda^{r}T^*_{\widehat\cX^a_k},$$
hence it does restrict to a section of
$K_{\widehat W_t}$ on $\smash{\widehat Z_t}$ when the latter is equipped
with its corres\-ponding multiplier ideal sheaf. A~priori, (5.1.16)
seems to be dependent on the choice of our liftings~$\widetilde\tau_j$, $u'$ and
$\xi$, but as it coincides with the intrinsically defined morphism 
$(5.1.14_k)$ at
the generic point of~$\cZ$, it must be uniquely defined. We further
pass to the integral closures on both sides of (5.1.16) to reach what 
we defined to be the multiplier ideals, according to~Prop.~1.2.
The only potential trouble to complete the proof of Lemma~5.1.13 is 
the presence of the divisor~$(\psi\circ\mu)^*\Delta$ in (5.1.17). This is overcome, 
at least in the special case we need, by the next lemma.\qed\vskip5pt

\claim 5.1.18. Lemma|Assume that $(\cZ,\cW)$ is induced by the $k$-th stage
$(\cZ_k,\cW_k)$ of a directed subvariety $(\cZ_0,\cW_0)\subset(\cX,T_{\cX/S})$
such that $\pi_{k,k-1}(\cZ)=\cZ_{k-1}$. If $N\ge(k+1)(n-1)$, one can arrange 
the choice of $\cZ^a\subset\cX^a_k$ in Lemma~{\rm 5.1.7} and of the
desingularization $\mu$ in Lemma~{\rm 5.1.13} in such a way that the
generic fiber $\smash{\widehat Z_t}$ does not meet the support of
$(\psi\circ\mu)^*\Delta$ 
in $(5.1.17)$ $($hence the projection of $(\psi\circ\mu)^{-1}(\Delta)\subset
\smash{\widehat\cZ\subset \widehat\cX^a_k}$ on $S$ will be contained
in a proper algebraic subvariety of~$S\,)$. We~then say that the
extension $\cZ\rightsquigarrow\cZ^a$ is equisingular.
\endclaim

\proof One can first take a desingularization $\mu_0:\widehat\cZ_0\to\cZ_0$ 
of~$\cZ_0$, and find a Zariski open set $S'\subset S$ over which 
$\pi\circ\mu_0:\widehat\cZ_0\to S$ is a submersion. We replace the family
$\cX\to S$ by its restriction
$\smash{\widehat\cZ'_0}:=(\pi\circ\mu_0)^{-1}(S')\to S'$
and observe that it is enough to extend $\cZ$ into $\cZ^a$ within
the relative and absolute Semple towers of 
$\smash{\widehat\cZ'_0}\to S'$ (an arbitrary extension over
$\smash{\widehat\cX^a_{k-1}\supset\widehat\cZ^a_{k-1}}$ will then do).
Therefore, it is enough to prove Lemma~5.1.18 when $\cZ_0=\cX$ and
$\pi_{k,k-1}(\cZ)=\cX_{k-1}$ (after a replacement of our ambient space $\cX$
by~$\cX'=\smash{\widehat\cZ'_0}$ of smaller dimension $n'\le n$). 
We set $r=\rank\cW$ and notice that 
$$\dim\cZ=\dim\cX_{k-1}+\rank T_{\cZ/\cX_{k-1}}=\dim\cX_{k-1}+r-1.$$ 
On the other hand, we have to take
$$\dim\cZ^a=\dim\cX^a_{k-1}+N+r-1<\dim\cX^a_k=\dim\cX^a_{k-1}+N+n-1,
$$
hence the codimension of $\cZ$ in $\cZ^a$ is extremely large, while 
the codimension of $\cZ^a$ in $\cX^a_k$ is small, equal to $n-r$ (this makes
the result not so surprising).

For the sake of explaining the argument in simpler terms, first assume 
that $\cZ$ and $\cW$
are nonsingular. We produce $\cZ^a\subset\cX^a_k$ as a complete intersection
$\{\sigma_j=0\}$ of codimension $n-r$ in $\cX^a_k$. We claim that we can 
take $\cZ^a$ to be nonsingular along any fixed smooth fiber $Z_t$, $t=t_0$,
if $N\ge(k+1)(n-1)$. In fact, if $\cI_{\cZ}$ is the ideal sheaf of $\cZ$ in
$\cX^a_k$, we can take the sections 
$\sigma_j\in H^0(\cX^a_k,\cI_\cZ\otimes A)$
in a sufficiently ample line bundle $A$ on~$\cX^a_k$, chosen so that global 
sections  $\sigma$ vanishing along $\cZ$ still generate all possible 
$1$-jets at any point of~$Z_t$, i.e.\ their differentials $d\sigma$ 
generate the conormal bundle~$\cN^*_{\cZ|Z_t}=(T_{\cX^a_k}/T_{\cZ})^*_{|Z_t}
=(\cI_{\cZ}/\cI_{\cZ}^2)_{|Z_t}\,$; it 
is enough to take $A\gg 0$ so that $(\cN^*_{\cZ}\otimes A)_{|Z_t}$ is generated 
by sections and the groups $H^1(\cX^a_k,\cI_{\cZ}^2\otimes A)$,
$H^1(\cX^a_k,\cI_{Z_t}\otimes\cN^*_\cZ\otimes A)$ vanish, thanks to the
exact sequences
$$
\eqalign{
&H^0(\cX^a_k,\cN^*_{\cZ}\otimes A)\to H^0(Z_t,(\cN^*_{\cZ}\otimes A)_{|Z_t})\to
H^1(\cX^a_k,\cI_{Z_t}\otimes\cN^*_\cZ\otimes A)=0,\cr
&H^0(\cX^a_k,\cI_{\cZ}\otimes A)\to H^0(\cX^a_k,\cI_{\cZ}/\cI_{\cZ}^2\otimes A)
\to H^1(\cX^a_k,\cI_{\cZ}^2\otimes A)=0.\cr}
$$
The conormal bundle $\cN^*_{\cZ}$ projects surjectively onto the vector
bundle $(\smash{\widetilde\cV}_k/\cW)^*$
of rank \hbox{$N+n-r$}. What we want is that the differentials 
$d\sigma_1,\,\ldots,\,d\sigma_r$ be pointwise linearly independent 
as sections of $(\smash{\widetilde\cV}_k/\cW)^*$ along 
$Z_t\,$; we then get transverse intersections
$$
\cW=\cV_k\cap\{d\sigma_j=0\},\qquad
\widetilde\cW=\widetilde\cV_k\cap\{d\sigma_j=0\}.
$$
Now, $\dim Z_t=n+(k-1)(n-1)+(r-1)\le (k+1)(n-1)$, and by a well 
known argument,
our condition on the independence of $n-r$ sections is true for generic
sections in a spanned vector bundle of 
$\rank\ge \dim Z_t+(n-r)$, i.e.\
$N\ge \dim Z_t$. Then $\smash{\widetilde\cW}$ will be
nonsingular, and $\smash{\widetilde\cW}/\cW\simeq\smash{\widetilde\cV}_k/
\cV_k=\cS_k$ along $Z_t$, as desired.

When $\cZ$ and $\cW$ are singular, we take an embedded desingularization 
$\mu:\smash{\widehat\cZ}\to\cZ\subset\cX^a_k$ in such a way that 
$\mu^*\cW$ becomes a nonsingular subbundle
of $\mu^*\cV_k$ (for this, we resolve the indeterminacies of the 
meromorphic map $\varphi$ to the relevant Grassmannian bundle, as 
already explained in Lemma~5.1.13). We work on a fixed nonsingular 
fiber $\smash{\widehat Z}_t$. What we want is that
the pull-backs $\mu^*d\sigma_j$, $1\le j\le n-r$ still cut out
a nonsingular rank $N$ vector subbundle in the rank $N+n-r$ bundle
$(\smash{\mu^*\widetilde\cV}_k/\mu^*\cW)^*$, in restriction to
the fiber $\smash{\widehat Z}_t$. This subbundle will
become our $\mu^*\smash{\widetilde\cW}$ [and as it will be already nonsingular,
there will be no  need to resolve further the indeterminacies of
$\mu^*\smash{\widetilde\cW}$ along the given fiber
$\smash{\widehat Z_t}$]. Moreover $(\psi\circ\mu)^{-1}(\Delta)$ will not meet 
$\smash{\widehat Z_t}$ by construction. The embedded resolution of
$\cZ$ in $\smash{\widehat\cX^a_k}$ yields 
$\mu^*\cI_{\cZ}=\cI_{\widehat\cZ}\cdot\cO(-\sum m_jE_j)$ where 
$\smash{\widehat\cZ}$ is the strict transform of $\cZ$ and the
$E_j$'s are the exceptional divisors, which we can assume to be normal
crossing and intersecting $\smash{\widehat\cZ}$ transversally.
The pull-back $\mu^*\sigma$ of any section $\sigma\in H^0(\cX^a_k,A\otimes
\cI_{\cZ})$ can be seen as a section $\widehat\sigma
\in H^0(\smash{\widehat\cX^a_k},\mu^*A\otimes\cO(-\sum m_jE_j))$
with coefficients in $\cI_{\widehat\cZ}$. By taking $A$ ample enough 
(and killing ad hoc direct image sheaves via~$\mu_*$), we
can achieve that the differentials $d\widehat\sigma$ still generate
along $\smash{\widehat Z}_t$ the normal 
bundle of $\smash{\widehat\cZ}$ in~$\smash{\widehat\cX^a_k}$. 
The same dimension argument as in the nonsingular case then allows us
to choose $\sigma_1,\,\ldots\,,\,\sigma_{n-r}$ so that
$\smash{\mu^*\widetilde\cW}/\mu^*\cW\simeq
\mu^*(\smash{\widetilde\cV}_k/\cV_k)=\mu^*\cS_k$ along $\smash{\widehat Z}_t$,
and we are done.\qed\vskip5pt

\noindent{\it End of proof of Theorem {\rm 5.1.}} If we combine 
(5.1.9) and Lemmas~5.1.13, 5.1.18 (thus,  under the assumption 
$N\ge(k+1)(n-1)$), we get a restriction morphism
$$
(\pi_{k,0}\circ\mu)^*\Lambda^{q'}T^*_{\cX|\widehat Z_t}
\longrightarrow K_{\widehat W_t}\otimes
\mu^*\big((\pi\circ\pi_{k,0})^*K_S\otimes\cO_{X_k}({\bf a}-N\,{\bf 1})
\big)_{|\widehat Z_t}.\leqno(5.1.19)
$$
Since $\Lambda^{q'}T^*_{\cX}$ is relatively ample over~$S$ and 
the factor $(\pi\circ\pi_{k,0})^*K_{S\,|Z_t}$ is trivial,
we see after multiplication by an additional small relatively ample term 
$\cO_{X_k}(\varepsilon{\bf c})$ that there is a weight ${\bf b}\in\bN^k$
such that $\smash{K_{\widehat W_t}\otimes\mu^*\cO_{X_k}({\bf b})_{|\widehat Z_t}}$
is big, hence $(Z_t,W_t)$ is of general type modulo
$X_{k,t}\to X_t$ for generic~$t$ (it is helpful here to know 
that $a_j>q-1\ge N$). This is a contradiction, hence $p_{k,\delta}$
is not dominant and $S_{k,\delta}=\smash{\overline{p_{k,\delta}(S)}}\subsetneq S$
(here we reproject down to $S$ in case there were a finite cover
$\smash{\widehat S^0}\to S$). If $N=\dim S$ is too small, we can artificially
increase the dimension of $S$ by replacing $S$ with $S\times\bC^m$,
$m\gg 1$, and put $X_{(t,t')}=X_t$, $\cZ_{(t,t')}=\cZ_t$, $(t,t')\in
S\times \bC^m$. The projection of the extended $p_{k,\delta}$'s to 
$S\times\bC^m$ is then non dominant, but since the property of 
$\cZ_{(t,t')}$ to be strongly of general type does not 
depend on~$t'$, we conclude that the projection to $S$ itself is
non dominant; in fact, the argument amounts to add many additional 
deformation parameters that are used to ``distort'' $\cZ^a$ along the 
fibers $\smash{\widehat Z_t}$, and so achieve the required equisingularity
property of $\cZ^a$ in Lemma~5.1.7. We conclude 
that $X_t$ is algebraically jet-hyperbolic for $t\in S\smallsetminus
\bigcup_{k,\delta}S_{k,\delta}$ and Theorem 5.1 is proved.\qed

\claim 5.2. Remarks|{\rm In fine, the main argument of the proof is the
existence of a non trivial morphism given by (5.1.9). If for all 
relevant subvarieties $\cZ\subset\cX_k$ one can find an ample subbundle
$\cA\subset\Lambda^{q'}T_{\cX}$ such that the composition
$$
(\pi_{k,0})^*\cA_{\,|\cZ}\longrightarrow
(\pi_{k,0})^*\Lambda^{q'}T^*_{\cX\,|\cZ}\longrightarrow K_{\widetilde\cW}
\otimes\cO_{\cX_k}({\bf a})_{|\cZ},\quad{\bf a}\in\bQ_+^k
$$
is non zero, then the conclusion still holds. This may allow to weaken 
the hypotheses on the positivity of $\Lambda^{q'}T_{\cX}$.

Also, one can in fact get the stronger conclusion that the very
general fiber $X_t$ is algebraically fully jet-hyperbolic in the sense
of Definition 4.1~(b). The only change is that we need a more sophisticated
version of Lemma~5.1.18; instead of just desingularizing
$\cZ_0\subset\cX$ and $\cZ\subset\cX_k$, we perform embedded 
desingularizations step by step for all intermediate projections 
$\cZ_\ell=\pi_{k,\ell}(\cZ)$, in case the ranks of
the associated directed structures $\cW_\ell$ gradually change. The arguments
are essentially the same and are left to the reader.\qed}
\endclaim

\noindent
{\bf 5.3. Universal family of complete intersections.} Let us consider
the universal family of complete intersections of dimension $n$, 
codimension $c$ and type $(d_1,\ldots,d_c)$ in complex projective 
$\bP^{n+c}=P(E)$, where $E\simeq\bC^{n+c+1}$ is a complex vector space. 
We can view it as a smooth family $\pi=\pr_1:\cX\to S$ 
where $S$ is a Zariski open set in 
$$\overline S=
\prod_{1\le j\le c} \Sym^{d_j}E^*\simeq\prod\bC^{N_j}=\bC^N,\qquad 
N_j={d_j+n+c\choose n+c},$$
and $\cX\subset S\times P(E)$ is the incidence variety defined by 
$$
\left\{
\eqalign{
&t=(t_1,\ldots,t_c)\in S,\quad z\in E\simeq\bC^{n+c+1},\quad 
t_j\simeq(t_{j,\alpha})\in\Sym^{d_j}E^*,\cr
\noalign{\vskip6pt}
&P_j(t,z):=t_j\cdot z^{d_j}=\sum_{|\alpha|=d_j}t_{j,\alpha}z^\alpha,\quad 1\le j\le c,
\quad\alpha=(\alpha_\ell)\in\bN^{n+c+1},\cr
\noalign{\vskip-3pt}
&\cX=\big\{(t,[z])\in S\times P(E)\,;\; P_j(t,z)=0,~1\le j\le c\big\}.\cr}
\right.\leqno(5.3.1)
$$ 
We denote by $\pr_1:\cX\to S$ and $\pr_2:\cX\to P(E)\simeq\bP^{n+c}$ the
natural projections. Here $S$ is the set of coefficients $t\in\bC^N$ that 
define a nonsingular subvariety~$X_t=\pr_1^{-1}(t)$ of codi\-mension~$c$
in $\bP^{n+c}$ (or rather $\{t\}\times\bP^{n+c}$). Notice that there is a 
natural action of $\GL(E)=\GL(n+c+1,\bC)$ on $\cX$ defined~by
$$
g\cdot((t_j),[z])=((t_j\circ g^{-1}),[g\cdot z]),
\quad g\in \GL(E),~~t_j\in\Sym^{d_j}E^*,\leqno(5.3.2)
$$
which simply consists of transforming the equations via an arbitrary linear
change of coordinates. We use the
following famous result proved by Claire Voisin [Voi96, Corollary~1.3]
(with the substitution of notation $k\mapsto c$, $n\mapsto n+c$, 
$l\mapsto N+n-q$ in our setting).

\claim 5.4. Proposition {\rm([Voi96])}|Over any affine Zariski open 
set $S^0\subset S$, the twisted tangent bundle 
$T_{\cX}\otimes\pr_2^*\cO_{\bP^{n+c}}(1)$ is generated by sections. Moreover,
the vector bundle $\Lambda^qT^*_{\cX}$ 
is generated by sections for $\sum d_j\ge 2n+c+N+1-q$,
and it is relatively very ample with respect to the projection $\cX\to S$
for $\sum d_j>2n+c+N+1-q$.
\endclaim

\proof Since the argument can be made very simple and very 
short, we give it here for the sake of completeness. If 
$(\varepsilon_j)_{0\le j\le n+c}$ denotes the
canonical basis of $\bZ^{n+c+1}$, we get sections of 
the tangent bundle of ${\rm cone}(\cX)$ over $\cX$ in 
$S\times E\simeq S\times\bC^{n+c+1}$
by taking the explicit vector fields
$$\leqalignno{
\kern40pt&\xi_{\ell,m}:=
z_m\frac{\partial~~\strut}{\partial z_\ell}-\sum_{1\le j\le c,\,|\alpha|=d_j}\alpha_\ell t_{j,\alpha}
\frac{\partial\kern40pt\strut}{\partial t_{j,\alpha-\varepsilon_\ell+\varepsilon_m}},
\quad 0\le\ell,m\le n+c,~~\alpha\in\bN^{n+c+1},&(5.4.1)\cr
&\eta_{j,\alpha,\ell,m}:=
z_m\frac{\partial\kern10pt\strut}{\partial t_{j,\alpha}}-
z_\ell\frac{\partial\kern40pt\strut}
{\partial t_{j,\alpha-\varepsilon_\ell+\varepsilon_m}},\quad 
|\alpha|=d_j,~~\alpha_\ell>0,~~0\le \ell\ne m\le n+c,&(5.4.2)\cr}
$$
which all yield zero derivative when applied to any of the polynomials
$P_j(t,z)$. In fact the vector fields (5.4.1) are just the Killing vector 
fields induced by the action of $\GL(E)$ on ${\rm cone}(\cX)$.
The natural $\bC^*$~action defined
by $\lambda\cdot (t,z)=(t,\lambda z)$ has an associated Euler vector field
$\varepsilon=\sum_{0\le \ell\le n+c}\,z_\ell\,\partial/\partial z_\ell$.
By taking the quotient with the rank $1$ subbundle $\cO_\cX\cdot\varepsilon$,
the $\xi_{\ell,m}$'s actually 
define sections of $T_{\cX}$ (homogeneity degree $0$ in $z$), while the
$\eta_{j,\alpha,\ell,m}$'s define sections of 
$T_{\cX}\otimes\pr_2^*\cO_{\bP^{n+c}}(1)$ (homogeneity degree $1$ in $z$).
We claim that the vector fields 
$$(t,[z])\mapsto\xi_{\ell,m}z_p~\mod\cO_\cX\cdot\varepsilon,\quad 
(t,[z])\mapsto \eta_{j,\alpha,\ell,m}~\mod\cO_{\cX}\cdot\varepsilon$$
generate $T_{\cX}\otimes\pr_2^*\cO_{\bP^{n+c}}(1)$ at every point.
In fact, as the $\xi_{\ell,m}$ already provide all ``vertical'' 
$z$-directions, we need only check that is is enough to add one 
(non tangent) ``horizontal'' vector field 
$\partial/\partial\smash{t_{j,\alpha_0^j}}$ for each $j=1,\ldots,c$ 
to generate the whole
ambient tangent space~$T_{\bP^{n+c}\times\bC^N,x}\,$, since the claim
then follows by a trivial (co)dimension argument. At a point
$(t,[z])$ where $z_0\ne 0$ (say), we take 
$\alpha_0^j=(d_j,0,\ldots,0)\in\bN^{n+c+1}$. Together with
$\partial/\partial\smash{t_{j,\alpha_0^j}}$, the vector fields
$$
z_0^{-1}\eta_{j,\alpha,\ell,0}:=
\frac{\partial\kern10pt\strut}{\partial t_{j,\alpha}}-
z_0^{-1}z_\ell\frac{\partial\kern40pt\strut}
{\partial t_{j,\alpha-\varepsilon_\ell+\varepsilon_0}},\quad 
\alpha_\ell>0,\quad\ell\ne 0\leqno(5.4.3)
$$
then generate $\smash{T_{\bC^{N_j}}}$ by a simple triangular 
matrix argument (increase the value of the \hbox{$0$-th} component of
$\alpha$ and decrease the $\alpha_\ell$'s, $\ell\ne 0$, until
$\alpha=\smash{\alpha_0^j}$). Let $\cO_S(-1)$ be the tautological line bundle
on $S$ (coming from the tautological line
bundle $\cO_{P(\overline S)}(-1)$ on $P(\overline S)$). Since  
$K_{S\times\bP^{n+c}}=\pr_1^*K_S\otimes\pr_2^*\cO_{\bP^{n+c}}(-n-c-1)$ and
$\cX$ is defined by 
sections of the line bundles 
$\pr_1^*\cO_S(-1)\otimes\pr_2^*\cO_{\bP^{n+c}}(d_j)$,
the adjunction formula gives
$$
K_\cX=\Lambda^{N+n}\,T^*_\cX=\cL_S\otimes\pr_2^*\cO_{\bP^{n+c}}
\big(\,{\textstyle\sum d_j}-n-c-1\big).
$$
where $\cL_S$ is the line bundle $\pr_1^*(\det T^*_S\otimes\cO_S(-c))$
(this bundle plays no role in the sequel since it can be made 
trivial by restricting $S$ to a suitable affine chart). Therefore
$$
\leqalignno{\strut\kern38pt
\Lambda^q\,T^*_\cX&=K_\cX\otimes\Lambda^{N+n-q}\,T_\cX&(5.4.4)\cr
&=\cL_S\otimes\pr_2^*\cO_{\bP^{n+c}}\big(\,{\textstyle\sum d_j}-n-c-1\big)
\otimes\Lambda^{N+n-q}\,T_\cX\cr
&=\cL_S\otimes
\pr_2^*\cO_{\bP^{n+c}}\big(\,{\textstyle\sum d_j}-2n-c-N-1+q\big)
\otimes\Lambda^{N+n-q}\,\big(T_\cX\otimes\pr_2^*\cO(1)\big).\cr}
$$
As $T_\cX\otimes\pr_2^*\cO(1)$ is generated by sections, Prop.~5.4
follows immediately.\qed
\vskip5pt

If we want the relative ampleness of $\Lambda^qT^*_{\cX}$ to hold
for $q>N$, we need $\sum d_j\ge 2n+c+1$. Theorem~5.1 then implies:

\claim 5.5. Corollary (solution of the Kobayashi conjecture)|For all
$n,c\ge 1$ and $d_j$ such that $\sum d_j\ge 2n+c+1$, the~very general complete
intersection of type $(d_1,\ldots,d_c)$ in complex projective space
$\bP^{n+c}$ is algebraically jet-hyperbolic, and thus Koba\-yashi hyperbolic. 
\endclaim

The simplest non trivial situation is the surface case $n=2$ in
codimension $c=1$. We then obtain the Kobayashi hyperbolicity of a very 
general surface $X\subset\bP^3$ of degree $d\ge 6$. The result seems
to be new even in this case, although Duval [Duv04] has shown by
elementary means the existence of a hyperbolic sextic (from this, it
already follows that there is a family of hyperbolic sextics over an open
set of parameters in Hausdorff topology). Geng Xu [Xu95] has shown that a 
very general quintic surface $X$ does not contain curves of genus $g\le 2$,
but as far as we know, this is not enough to conclude that $X$ is
Kobayashi hyperbolic.

\claim 5.6. Remark|{\rm It would be good to know if Kobayashi hyperbolicity
is a Zariski open condition, in particular, whether one can replace 
``very general'' by ``general'' in Cor.~5.5. This would require further
investigations, but such a result might be accessible by taking into account
Remark~3.3, which shows that the ``bad sets'' $Z$ to consider are 
somehow bounded.}
\endclaim

\claim 5.7. Remark|{\rm In the case $n\ge 2$ and 
$\sum d_j=2n+c$ (and especially in the ``border case'' $d=2n+1$ 
of hypersurfaces), it follows from
Prop.~5.4 due to [Voi96] that $(\Lambda^qT^*_{\cX})_{|X_t}$ is generated by 
sections for $q=N+1$ and very ample for $q\ge N+2$. It would then be 
natural to look at the degeneration sets occurring for all
appropriate subvarieties $\cZ$ in the various stages 
of the relative Semple tower. The arguments used by Claire Voisin
([Voi96], Theorem~1.6 and its proof) indicate a possibility to
analyze the situation, but certain remaining degeneracies seem 
to require intricate Wronskian and flag manifold arguments.}
\endclaim

\claim 5.8. Case of complements|{\rm Our techniques also apply to
study the Kobayashi hyperbo\-li\-city of complements
$\bP^n\smallsetminus X$, when $X$ is an algebraic hypersurface of
degree~$d$ in $\bP^n$. In fact, if $X=\{P(z_0,\ldots,z_n)=0\}$, one
can introduce the hypersurface
$$Y=\{z_{n+1}^d-P(z_0,\ldots,z_n)=0\}\subset\bP^{n+1}.$$
It is trivial to show that the Kobayashi hyperbolicity of $X$ implies
the Kobayashi hyperbo\-licity of~$X$, since the natural projection
$$
\rho:Y\to\bP^n,\quad(z_0,\ldots,z_{n+1})\mapsto (z_0,\ldots,z_n)
$$
defines an unramified $d:1$ cover from $Y\smallsetminus\rho^{-1}(X)$ onto
$\bP^n\smallsetminus X$. We have a universal family $\cY\to S$ by
looking at the parameter space given by coefficients of~$P$. This is
just a subfamily of the universal family of degree $d$ hypersurfaces,
and we only have to check that Prop.~5.4 still applies when we have no 
dependence on the variable $z_{n+1}$ except for the monomial $z_{n+1}^d$.
Here the group acting on the ambient projective space $\bP^{n+1}$ 
is taken to be
$$
\GL(n+1,\bC)\times\bC^*\subset\GL(n+2,\bC),
$$
and one can see that the last Killing vector field
$z_{n+1}\partial/\partial z_{n+1}+(\ldots)$ introduces some degenerations
on $z_{n+1}=0$ -- and only there. We easily conclude by our techniques
that $\ECL(X)\subset X\cap\{z_{n+1}=0\}$ for $P$ very general of degree 
$d\ge 2n+2$, but since we also have $\ECL(Y)=\emptyset$, we conclude
that $\bP^n\smallsetminus X$ is Kobayashi hyperbolic for $X$ very general
of degree~$d\ge 2n+2$. Zaidenberg [Zai87] has shown that this conclusion fails
for $d=2n$. One could hope to improve the bound to $d\ge 2n+1$ by 
introducing logarithmic Semple jet bundles, as suggested by Dethloff and Lu
[DL01], and apply the idea suggested in Remark~5.7.}
\endclaim
\vskip8mm

\plainsection{References}

{\ninepoint
\parskip 1.5pt plus 0.5pt minus 0.5pt

\Bibitem[Cle86]&Clemens, H.:& Curves in generic hypersurfaces.&
Ann.\ Sci.\ \'Ecole Norm.\ Sup.\ {\bf 19} (1986), 629--636&

\Bibitem[Dem95]&Demailly, J.-P.:& Algebraic criteria for Kobayashi
hyperbolic projective varieties and jet differentials.& AMS Summer
School on Algebraic Geometry, Santa Cruz 1995, Proc.\ Symposia in
Pure Math., ed.\ by J.~Koll\'ar and R.~Lazarsfeld, Amer.\ Math.\ Soc., 
Providence, RI (1997), 285–-360&

\Bibitem[Dem97]&Demailly, J.-P.:& Vari\'et\'es hyperboliques et
\'equations diff\'erentielles alg\'e\-bri\-ques.& Gaz.\ Math.\ {\bf 73}
(juillet 1997) 3--23, and {\ninett \HOMEPAGE/cartan\_{}augm.pdf}&

\Bibitem[DEG00]&Demailly, J.-P., El Goul, J.:& Hyperbolicity of
generic surfaces of high degree in projective 3-space.& Amer.\ J.\ 
Math.\ {\bf 122} (2000) 515--546&

\Bibitem[Dem11]&Demailly, J.-P.:& 
Holomorphic Morse inequalities and the Green-Griffiths-Lang
conjecture.& November 2010, arXiv: math.AG/1011.3636, dedicated to the
memory of Ec\-kart Viehweg; Pure and Applied Mathematics 
Quarterly {\bf 7} (2011) 1165--1208&

\Bibitem[DL01]&Dethloff, G., Lu, S.S.Y.:& Logarithmic jet bundles and 
applications.& Osaka J.~Math.\ {\bf 38} (2001) 185--237&

\Bibitem[DMR10]&Diverio, S., Merker, J., Rousseau, E.:& Effective
algebraic degeneracy.& Invent.\ Math.\ {\bf 180} (2010) 161--223&

\Bibitem[DR13]&Diverio, S., and Rousseau, E.:& The exceptional set and 
the Green-Griffiths locus do not always coincide.& arXiv: math.AG/1302.4756 
(v2)&

\Bibitem[Duv04]&Duval, J.:& Une sextique hyperbolique dans $\hbox{\ninebb P}^3(\hbox{\ninebb C})$.& Math.\ Ann.\ {\bf 330} (2004) 473--476&

\Bibitem[Ein88]&Ein, L.:& Subvarieties of generic complete intersections.&
Invent.\ Math.\ {\bf 94} (1988), 163--169&

\Bibitem[Ein91]&Ein, L.:& Subvarieties of generic complete intersections, II.&
Math.\ Ann.\ {\bf 289} (1991), 465--471&

\Bibitem[GG79]&Green, M., Griffiths, P.:& Two applications of algebraic
geometry to entire holomorphic mappings.& The Chern Symposium 1979,
Proc.\ Internal.\ Sympos.\ Berkeley, CA, 1979, Springer-Verlag, New York
(1980), 41--74&

\Bibitem[Kob70]&Kobayashi, S.& Hyperbolic manifolds and holomorphic
mappings.& Volume 2 of Pure and Applied Mathematics. Marcel Dekker Inc., 
New York, 1970&

\Bibitem[Kob78]&Kobayashi, S.& Hyperbolic complex spaces.& Volume 318
of Grundlehren der Mathematischen Wissenschaften, Springer-Verlag, 
Berlin, 1998&

\Bibitem[Lan86]&Lang, S.:& Hyperbolic and Diophantine analysis.&
Bull.\ Amer.\ Math.\ Soc.\ {\bf 14} (1986) 159--205&

\Bibitem[McQ98]&McQuillan, M.:& Diophantine approximation and foliations.&
Inst.\ Hautes \'Etudes Sci.\ Publ.\ Math.\ {\bf 87} (1998) 121--174&

\Bibitem[McQ99]&McQuillan, M.:& Holomorphic curves on hyperplane sections 
of $3$-folds.& Geom.\ Funct.\ Anal.\ {\bf 9} (1999) 370--392&

\Bibitem[Pau08]&P\u{a}un, M.:& Vector fields on the total space of 
hypersurfaces in the projective space and hyperbolicity.& Math.\ Ann.\
{\bf 340} (2008) 875--892&

\Bibitem[Siu02]&Siu, Y.T.:& Some recent transcendental techniques in
algebraic and complex geometry.& In: Proceedings of the International
Congress of Mathematicians, Vol.~I (Beijing, 2002), Higher Ed.\ Press, Beijing,
2002, 439--448&

\Bibitem[Siu04]&Siu, Y.T.:& Hyperbolicity in complex geometry.& In: The 
legacy of Niels Henrik Abel, Springer, Berlin, 2004, 543--566&

\Bibitem[Siu12]&Siu, Y.T:& Hyperbolicity of generic high-degree hypersurfaces
in complex projective spaces.& arXiv: math.CV, 1209.2723v3, 87 pages&

\Bibitem[SY96]&Siu, Y.T., Yeung, S.K.:& Hyperbolicity of the complement of
a generic smooth curve of high degree in the complex projective plane.& Invent.\
Math.\ {\bf 124} (1996), 573--618&

\Bibitem[Voi96]&Voisin, C.:& On a conjecture of Clemens on rational curves
on hypersurfaces.& J.~Differential Geometry {\bf 44} (1996) 200--213&

\Bibitem[Voi98]&Voisin, C.:& A correction on ``A conjecture of Clemens 
on rational curves on hypersurfaces''.& J.~Differential Geometry {\bf 49} 
(1998) 601--611&

\Bibitem[Xu94]&Xu Geng:& Subvarieties of general hypersurfaces in projective
space.& J.~Differential Geom.\ {\bf 39} (1) (1994) 139--172&

\Bibitem[Zai87]&Zaidenberg, M.:& The complement to a general hypersurface of degree $2n$ in $\hbox{\ninebb CP}^n$ is not hyperbolic.& (Russian) Sibirsk.\ Mat.\ Zh.\
{\bf 28} (3) (1987) 91--100&

}

\vskip3mm\noindent
Jean-Pierre Demailly\\
Institut Fourier, Universit\'e Grenoble-Alpes\\
BP74, 100 rue des Maths, 38402 Saint-Martin d'H\`eres, France\\
\emph{e-mail}\/: jean-pierre.demailly@ujf-grenoble.fr

\end{document}